\documentclass{amsart}
\usepackage[utf8]{inputenc}

\usepackage{amsthm}
\usepackage{amsfonts}
\usepackage{mathrsfs}
\usepackage{graphicx}

\usepackage{xcolor}
\usepackage{hyperref}

\usepackage{amssymb}
\usepackage{faktor}

\usepackage{caption}
\usepackage{subcaption}

\title{Piecewise circular curves and Positivity}
\author{Jean-Philippe Burelle}
\thanks{This project has received funding from the European Research Council (ERC) under the European Union’s Horizon 2020 research and innovation programme (ERC starting grant DiGGeS, grant agreement No 715982). The first author acknowledges the support of the Natural Sciences and Engineering Research Council of Canada (NSERC), [funding reference number RGPIN-2020-05557]
}
\author{Ryan Kirk}

\newtheorem{lem}{Lemma}[section]
\newtheorem{prop}[lem]{Proposition}
\newtheorem{cor}[lem]{Corollary}
\newtheorem{thm}[lem]{Theorem}
\newtheorem{example}[lem]{Example}

\theoremstyle{remark}
\newtheorem{rmk}[lem]{Remark}

\theoremstyle{remark}

\theoremstyle{definition}
\newtheorem{defn}[lem]{Definition}

\renewcommand{\vec}[1]{\mathbf{#1}}

\newcommand{\V}{\mathsf{V}}
\newcommand{\w}{\omega}
\newcommand{\Proj}{\mathbb{P}}

\newcommand{\bR}{\mathbb{R}}
\newcommand{\bC}{\mathbb{C}}
\newcommand{\RP}{\mathbb{RP}}
\newcommand{\CP}{\mathbb{CP}}
\newcommand{\Span}{\mathrm{span}}

\newcommand{\Flag}{\mathsf{Flag}}

\newcommand{\SO}{\mathsf{SO}}
\newcommand{\PSp}{\mathsf{PSp}}
\newcommand{\Sp}{\mathsf{Sp}}

\newcommand{\PSL}{\mathsf{PSL}}
\newcommand{\M}{\mathcal{M}}
\newcommand{\Lag}{\mathsf{Lag}}

\newcommand{\diag}{\mathrm{diag}}

\newcommand{\Hom}{\mathrm{Hom}}

\newcommand{\dif}{\mathrm{d}}

\begin{document}
\maketitle

\begin{abstract}
    We introduce the moduli space of generic piecewise circular $n$-gons in the Riemann sphere and relate it to a moduli space of Legendrian polygons. We prove that when $n=2k$, this moduli space contains a connected component homeomorphic to the Fock-Goncharov space of $k$-tuples of positive flags for $\PSp(4,\bR)$ and hence is a topological ball. We characterize this component geometrically as the space of simple piecewise circular curves with decreasing curvature.
\end{abstract}

\section{Introduction}

A \emph{piecewise circular curve} is a curve in the plane made of finitely many arcs of circles, such that tangents agree at the intersection of pieces. We will call a closed piecewise circular curve made of $n$ circular segments a circular $n$-gon. Since Möbius transformations map circular arcs to circular arcs or line segments, it makes sense to consider piecewise circular curves in the Riemann sphere up to Möbius transformations, allowing circular and linear pieces.

A co-orientation on a piecewise circular curve is a continuous choice of perpendicular orientation on the curve. There is a natural map on the set of co-oriented piecewise circular curves which moves each point on the curve a fixed distance $d$ in the direction of the co-orientation. We will call this map the radial translation with parameter $d$. If the piecewise circular curve is interpreted as a wavefront, this map models wave propagation with $d$ being the time parameter.

We consider the natural problem of classifying co-oriented piecewise circular curves up to Möbius transformations and radial translations.

Piecewise circular curves have been previously been investigated, especially in the setting of curve approximation (\cite{BanchoffGiblin1993}, \cite{BanchoffGiblin1994}, \cite{rosin1999}) but as far as we know this classification problem has not been considered in the literature.

Since the radial translation may introduce singularities by collapsing a circular arc to a point, we allow such singularities in our piecewise circular curves. Keeping track of the co-orientation, collapsed circular arcs are not degenerate as curves in the space of oriented contact elements of the sphere, and the action of radial translations and Möbius transformation is well-defined.

\begin{figure}[h]
    \centering
    \begin{subfigure}[b]{0.3\textwidth}
         \centering
         \includegraphics[width=\textwidth]{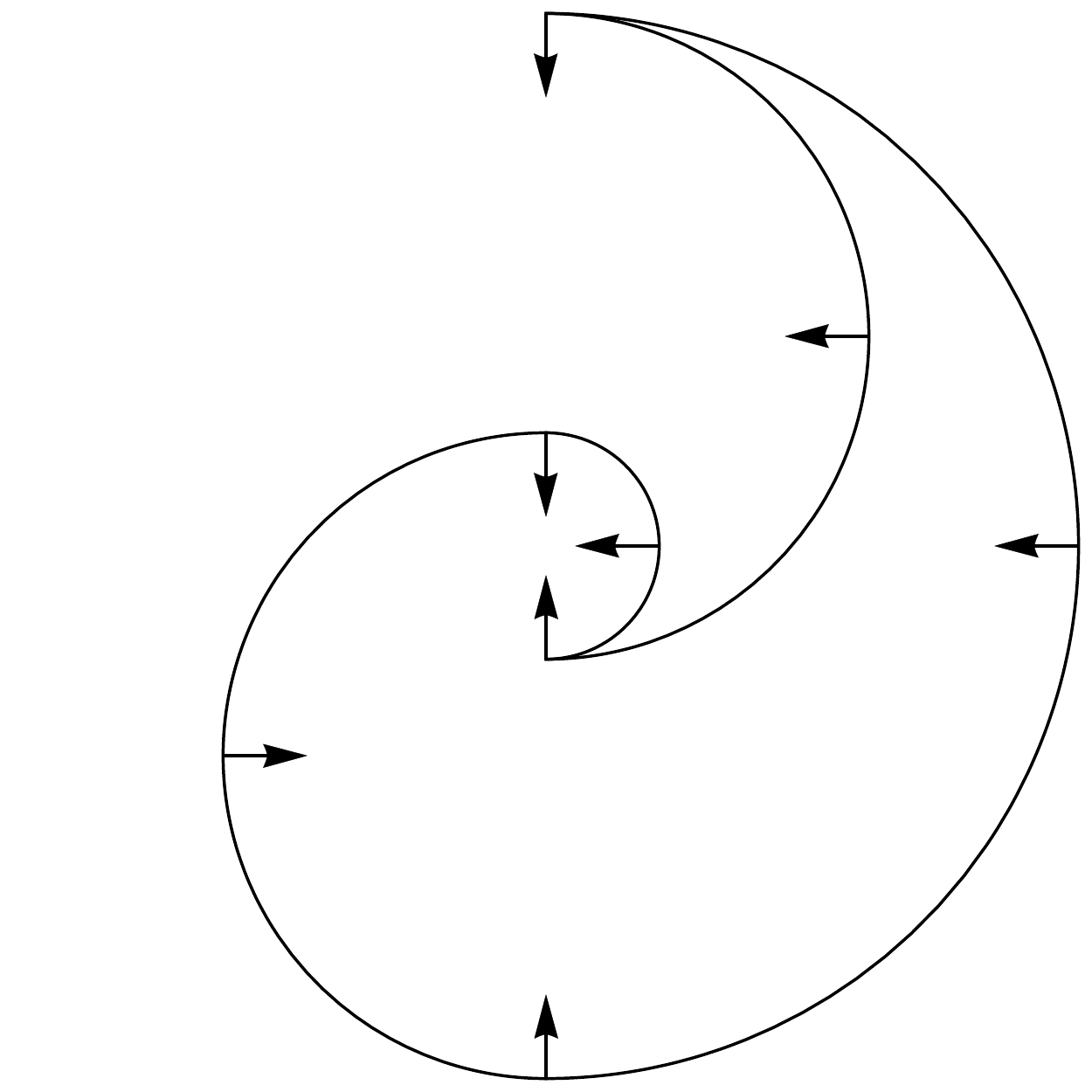}
     \end{subfigure}
     \hfill
     \begin{subfigure}[b]{0.3\textwidth}
         \centering
         \includegraphics[width=\textwidth]{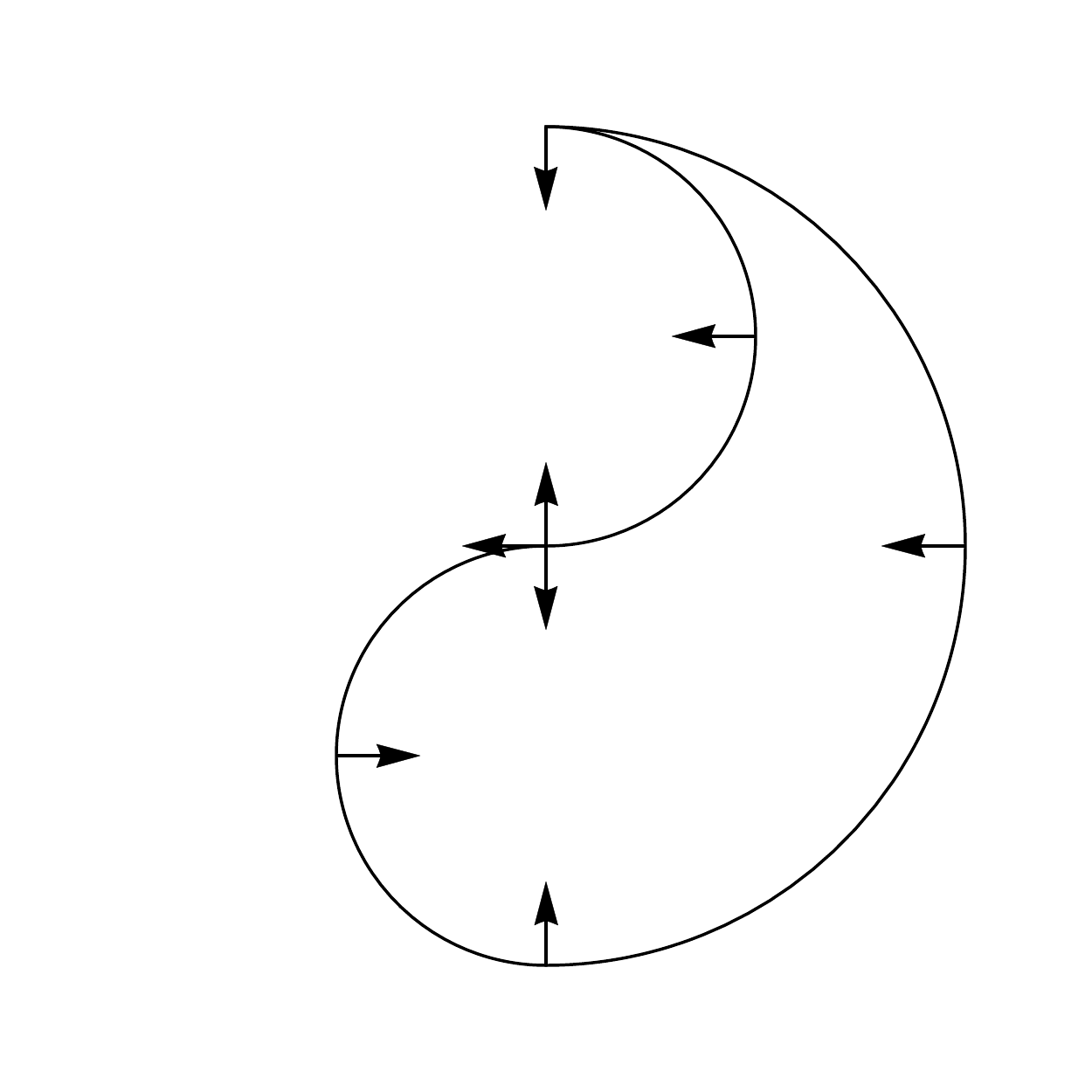}
     \end{subfigure}
     \hfill
     \begin{subfigure}[b]{0.3\textwidth}
         \centering
         \includegraphics[width=\textwidth]{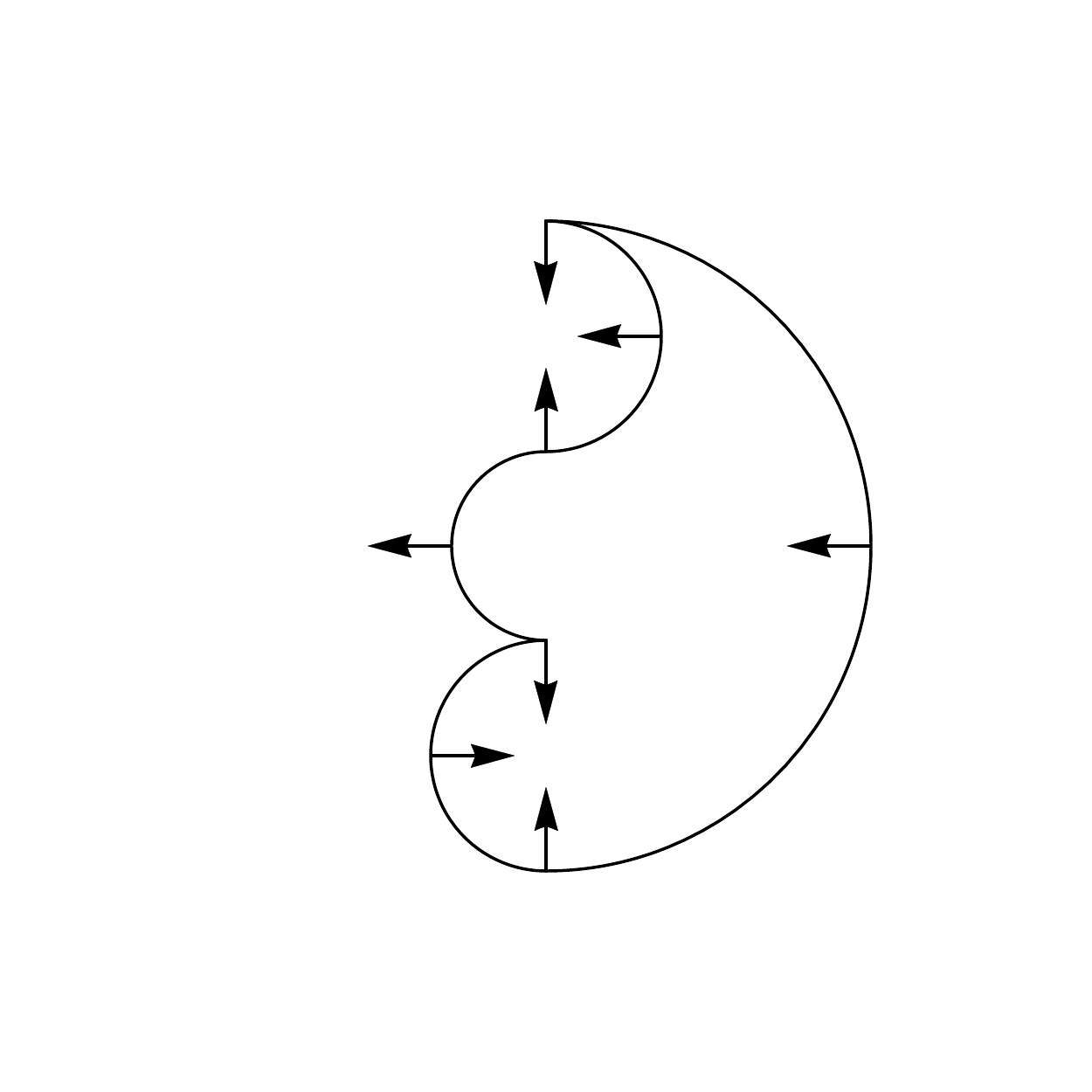}
     \end{subfigure}
        \caption{A simple co-oriented piecewise circular quadrilateral and its images under radial translation. The central figure shows a collapsed circular arc.}
        \label{fig:threegraphs}
\end{figure}
Even though we cannot compare the radii of any two co-oriented circles in the Riemann sphere, there is a partial cyclic order coming from the Maslov index which detects when three co-oriented circles are nested, and is invariant under Möbius transformations and radial translations. When the collection of circles extending the edges of a co-oriented piecewise circular curve is cyclically ordered, we will say that this curve has decreasing curvature. This property means that piecewise circular curves with decreasing curvature are locally spirals.

We say that a co-oriented piecewise circular curve is \emph{generic} if non-adjacent vertices are never part of a common co-oriented circle, and non-adjacent edges are never arcs of tangent co-oriented circles. We say it is \emph{simple} if \emph{all} of its radial translates have no self-intersections. This is equivalent to there being no co-oriented circle which is tangent to two non-adjacent edges (see Prop. \ref{prop:simpleifftransverse}).

We define the moduli spaces rigorously and make the following simple observations.
\begin{prop}
  The moduli space of generic circular $n$-gons with $n\ge 5$ is a smooth $2(n-5)$-dimensional manifold.
\end{prop}

\begin{prop}
  The moduli space of simple circular $4$-gons consists of two points.
\end{prop}

\begin{prop}
  As a curve in the space of oriented contact elements to the $2$-sphere, a simple piecewise circular polygon is non-contractible.
\end{prop}

A related moduli space, the space of \emph{$(2,n)$-Lagrangians configurations} in $\bR^4$ was investigated in \cite{ConleyOvsienko} and \cite{MorierGenoud1}. The precise relation is:

\begin{prop}
  The moduli space of simple circular $n$-gons is a two-to-one trivial cover of the space $\mathcal{L}(2,n)$ of generic configurations of $n$ Lagrangians in $\bR^4$.
\end{prop}

Among generic circular $2k$-gons, there is a particular class consisting of curves which are simple and have decreasing curvature. Our main result is:

\begin{thm}
  Let $k\ge 1$. In the moduli space of generic circular $2k$-gons, the subspace of simple curves with decreasing curvature is a connected component homeomorphic to an open ball.
\end{thm}

We also show that for hexagons, the curvature condition is not necessary, a result which was proved in the second author's Ph.D thesis \cite{kirkthesis}:
\begin{prop}
  A simple circular hexagon either has increasing or decreasing curvature.
\end{prop}

The strategy of the proof is to construct a homeomorphism between the moduli space of pairwise transverse $k$-tuples of $\PSp(4,\bR)$ flags and the moduli space of generic circular $2k$-gons. Under this homeomorphism, positive $k$-tuples of flags are in bijection with positive circular $2k$-gons. Fock and Goncharov showed in \cite{FG2006} that the subspace of positive $k$-tuples of flags is a connected component homeomorphic to a ball, and so the theorem follows.

Positive $k$-tuples of flags for the group $\PSL(3,\bR)$ also parameterize a natural moduli space of geometric objects: pairs of nested convex polygons in the projective plane \cite{FGprojective}. Our result provides a similar interpretation for positive $k$-tuples of flags in $\PSp(4,\bR)$. Both can be considered toy models of higher Teichm\"uller spaces.

As a consequence of the techniques used to prove the main theorem, we also obtain the following result about positive curves and Hitchin representations of surface groups into $\PSp(4,\bR)$:

\begin{thm}
  Positive Frenet $\PSp(4,\bR)$ flag curves are lifts of simple closed curves with decreasing curvature on the sphere.
\end{thm}

\begin{cor}
  Let $\Sigma$ be a closed surface and $\rho : \pi_1(\Sigma) \rightarrow \PSp(4,\bR)$ a Hitchin representation. Then, the limit curve of $\rho$ is the contact lift of a simple closed curve with decreasing curvature on the sphere.
\end{cor}

This corollary gives an alternative geometric interpretation of the $\PSp(4,\bR)$ Hitchin component to the one presented in \cite{guichardwienhard}. Guichard and Wienhard interpret the Hitchin component as the moduli space of convex-foliated contact projective structures on the tangent bundle of $\Sigma$. Corollary 1.8 gives an interpretation of this component as the moduli space of $\pi_1(\Sigma)$-equivariant simple closed curves with decreasing curvature on the $2$-sphere.

We now describe the structure of the paper. After establishing some notation, in Section \ref{sec:LegendrianPolygons} we introduce Legendrian polygons and their moduli spaces. In Section \ref{sec:circles} we develop the dictionary between the contact geometry of projective $3$-space and the geometry of circles on the $2$-sphere. Most of these results are known and have been described using the perspective of the Lie group $\SO^0(3,2)$ (for instance in \cite{cecil}). For our purposes, the group $\PSp(4,\bR)$ which is isomorphic to $\SO^0(3,2)$ is more convenient and so we develop \emph{Lie circle geometry} from that point of view in a coordinate-free way. Finally, in Section \ref{sec:Positivity} we recall the notion of positivity of a configuration of flags and prove the main theorem and describe the applications to higher Teichmüller theory.

\section{Notation}
When working with indexed objects $A_i$ with a natural cyclic structure (for instance vertices of a polygon) we will consider indices from the set $\{1,2,\dots,n\}$ \emph{modulo $n$}, meaning for instance that $A_{n+1}=A_1$ and $A_{n+2}=2$. When writing inequalities of the form $i<j$ for those indices, we use the unique representatives in $\{1,2,\dots,n\}$.

If $\V$ is a vector space over a field $k$, we denote by $\Proj(\V)$ the associated projective space. If $\mathsf{U}\subset \V$ is a vector subspace, we also denote by $\Proj(\mathsf{U})\subset \Proj(\V)$ the corresponding projective subspace. Over $\bR$, we also define the \emph{sphere of directions} $\mathbb{S}(\V) := (\V-\{0\})/\bR^+$ and we use the same notation $\mathbb{S}(E)$ for the sphere bundle associated to a real vector bundle $E\rightarrow M$ over a manifold $M$. If $\V$ is a vector space over two different fields, we will use $\Proj_k(\V)$ to specify which field we are considering. We use bold variables $\vec{v}\in\V$ to denote elements of $\V$ and brackets $[\vec{v}]_F \in \Proj_F(\V)$ to denote projective equivalence classes. 

A \emph{symplectic form} $\w$ on $\V$ is a skew-symmetric bilinear form. When $\V$ is equipped with a symplectic form, we will call it a symplectic vector space. Let $S\subset \V$ be a subset. Then,
\[S^\perp = \{\vec{v}\in \V ~|~ \w(\vec{v},\vec{u}) = 0,~ \forall \vec{u}\in S\}\]
denotes the orthogonal subspace to $S$.

A Lagrangian subspace $L\subset V$ is a subspace such that $L^\perp = L$.
The Lagrangian Grassmannian $\Lag(\V)$ is the space of Lagrangian subspaces in $\V$. When a basis of $\V$ is fixed and $\dim(\V)=2n$, we will use $2n\times n$ matrices to denote elements of $\Lag(\V)$, understood as maps $\bR^n\rightarrow \V$. This notation is unique up to right-multiplication by an invertible $n\times n$ matrix.

\section{Legendrian Polygons}\label{sec:LegendrianPolygons}
  Let $(\V,\w)$ be a $4$-dimensional symplectic vector space over $\bR$. A \emph{Lagrangian} in $\V$ is a $2$-dimensional subspace on which the symplectic form $\w$ vanishes.
  
  Orthogonality with respect to $\w$ defines a hyperplane distribution $\vec{v}^\perp \subset T_\vec{v} \V$. This distribution projects to a contact structure on the projective space $\Proj(\V)$, making it a contact manifold.
  
  \begin{defn}
    Two points $p,q \in \Proj(\V)$ will be called \emph{incident} whenever $p \in \Proj(q^\perp)$, or equivalently $q \in \Proj(p^\perp)$.
  \end{defn}
  
  A path in $\Proj(\V)$ which is always tangent to the contact distribution is called \emph{Legendrian}. 
  
  \begin{defn}
    A \emph{Legendrian polygon} in $\Proj(\V)$ is the image of a closed, piecewise linear path such that each segment is Legendrian, with labeled vertices.
  \end{defn}

  We will assume that all Legendrian polygons are nondegenerate in the following sense: no two adjacent segments are part of the same projective line and no segment is degenerated to a point.
  
  Let $(\vec{v}_1,\dots,\vec{v}_n)$ be nonzero vectors in $\V$ such that $\vec{v}_k,\vec{v}_{k+1},\vec{v}_{k+2}$ are linearly independent, and $\w(\vec{v}_k,\vec{v}_{k+1})=0$ for $1\le k\le n$. We define two Legendrian polygons $P_\pm(\vec{v}_1,\dots,\vec{v}_n)$ with $n$ segments $\gamma_1,\dots,\gamma_n$ given by:
  \[\gamma_k(t) = [(1-t)\vec{v}_k + t \vec{v}_{k+1}]\]
  for $0\leq t\leq 1$ and $1\le k \le n-1$ and
  \[\gamma_n(t) = [(1-t)\vec{v}_n \pm t \vec{v}_1].\]
  Any Legendrian polygon can be parametrized this way, and the ambiguity in the choice of representatives $\vec{v}_1,\dots,\vec{v}_n$ is multiplication of each $\vec{v}_j$ by a positive scalar $\lambda_j>0$ and simultaneously multiplying all $\vec{v}_j$ by $-1$.
  
  \begin{prop}\label{prop:pmcontractible}
    For any nonzero vectors $\vec{v}_1,\dots, \vec{v}_n\in \V$, the polygon $P_+(\vec{v}_1,\dots,\vec{v}_n)$ is contractible, and
    the polygon $P_-(\vec{v}_1,\dots,\vec{v}_n)$ generates $\pi_1(\Proj(\V))$.
    \begin{proof}
    
      Let $[\vec{u}_1],\dots,[\vec{u}_n]$ be $n$ distinct points on a projective line $\ell\subset \Proj(\V)$ disjoint from the polygon, placed in that order. Multiplying certain $\vec{u}_k$ by $-1$ if needed, we may assume that the full line is parametrized by $P_-(\vec{u}_1,\dots,\vec{u}_n)$.
      
      For each vector $\vec{v}_i$, choose a path $\vec{v}_i(s)\in \V - \{0\}$ such that $\vec{v}_i(0)=\vec{v}_i$ and $\vec{v}_i(1) = \vec{u}_i$. Perturbing the path if needed, $\vec{v}_k(s)$ and $\vec{v}_{k+1}(s)$ are linearly independent for all $s$. Then, the homotopy between $P_\pm(\vec{v}_1,\dots,\vec{v}_n)$ and $P_\pm(\vec{u}_1,\dots,\vec{u}_n)$
      
      \[\gamma_k(t,s) := [(1-t)\vec{v}_k(s) + t \vec{v}_{k+1}(s)]\]
      \[\gamma_n(t,s) := [(1-t)\vec{v}_n(s) \pm t \vec{v}_{k+1}(s)]\]
      
      is well-defined. The path $P_+(\vec{u}_1,\dots,\vec{u}_n)$ is contained in a line segment, hence contractible, and the path $P_-(\vec{u}_1,\dots,\vec{u}_n)$ is a full projective line, hence generates $\pi_1(\Proj(\V))$.

      
    \end{proof}
  \end{prop}
  

  \begin{defn}
     A Legendrian polygon $P$ is \emph{generic} if its non-adjacent vertices are non-incident and its non-adjacent edges are parts of transverse projective lines.
     
     The moduli space of generic Legendrian $n$-gons up to the action of $\PSp(\V)$ will be denoted by $\mathscr{P}_n$.
     
     As a consequence of Proposition \ref{prop:pmcontractible}, $\mathscr{P}_n$ naturally separates into a disjoint union of the space of generic contractible $n$-gons $\mathscr{P}_n^+$ and the space of generic non-contractible $n$-gons $\mathscr{P}_n^-$, and each subspace $\mathscr{P}_n^+, \mathscr{P}_n^-$ is a union of connected components.
  \end{defn}
  
  \begin{defn}
  A Legendrian polygon $P$ is \emph{transverse} if no pair of points $p,q\in P$ belonging to distinct, non-adjacent closed segments is incident.
\end{defn}

We first observe that there are no Legendrian triangles, so the first interesting case of Legendrian polygon is $n=4$. This is because the vertices of a Legendrian triangle would span a $3$-dimensional totally isotropic subspace, which is impossible by nondegeneracy of $\w$.

As a first example, in a basis $\vec{e}_1,\vec{e}_2,\vec{e}_3,\vec{e}_4$ in which the symplectic form $\w$ is represented by the matrix
\begin{equation}\label{eqn:OmegaMatrix}
\Omega=\begin{pmatrix} 0 & 0 & 0 & 1\\0 & 0 & -1 & 0\\0 & 1 & 0 & 0\\-1 & 0 & 0 & 0\end{pmatrix},
\end{equation}
the Legendrian quadrilateral $P_-(\vec{e}_1,\vec{e}_2,\vec{e}_4,-\vec{e}_3)$ is generic and transverse. This is a consequence of the following lemma:

\begin{lem}\label{lem:positive_segment_pair}
  Let $\vec{u}_1,\vec{u}_2,\vec{v}_1,\vec{v}_2 \in \V$ such that $\w(\vec{u}_1,\vec{u}_2)=\w(\vec{v}_1,\vec{v}_2)=0$ and $\w(\vec{u}_1,\vec{v}_2)\neq 0$. Then, the points $[(1-t)\vec{u}_1 + t \vec{u}_2]$ and $[(1-s)\vec{v}_1 + s \vec{v}_2]$ are non-incident for all $t,s\in(0,1)$ if and only if the four symplectic products
  $\w(\vec{u}_1,\vec{v}_1)$, $\w(\vec{u}_1,\vec{v}_2)$, $\w(\vec{u}_2,\vec{v}_1)$, and $\w(\vec{u}_2,\vec{v}_2)$ are all nonpositive or all nonnegative.
  \begin{proof}
     The points $[(1-t)\vec{u}_1 + t \vec{u}_2]$ and $[(1-s)\vec{v}_1 + s \vec{v}_2]$ are incident if and only if $\w((1-t)\vec{u}_1 + t \vec{u}_2,(1-s)\vec{v}_1 + s \vec{v}_2) = 0$, which expands to
     \begin{equation*}
         (1-t)(1-s)\,\w(\vec{u}_1,\vec{v}_1) + (1-t)s~\,\w(\vec{u}_1,\vec{v}_2) + t(1-s)\,\w(\vec{u}_2,\vec{v}_1) + t s\,\w(\vec{u}_2,\vec{v}_2) = 0.
     \end{equation*}
     
     Suppose that the four symplectic products $\w(\vec{u}_i,\vec{v}_j)$ are all nonpositive or all nonnegative. Then, since one of them is nonzero the above equation has no solutions for $0 < t,s < 1$, and so $[(1-t)\vec{u}_1 + t \vec{u}_2]$ and $[(1-s)\vec{v}_1 + s \vec{v}_2]$ are non-incident.
     
     Conversely, note that for fixed values of $a,b,c,d$ the function
     \[\phi(t,s) = (1-t)(1-s)a + (1-t)s\,b + t(1-s)c + ts\,d\]
     satisfies $\phi(0,0)=a$, $\phi(0,1)=b$, $\phi(1,0)=c$, and $\phi(1,1)=d$.
     
     If any two of $a,b,c,d$ have different signs we can, by applying the intermediate value theorem to a path joining the vertices of the square $[0,1]\times[0,1]$, find $(t,s)\in (0,1)\times(0,1)$ such that $\phi(t,s)=0$.
  \end{proof}
\end{lem}

\begin{prop}\label{prop:pplusnottransverse}
  The Legendrian polygon $P_+(\vec{v}_1,\dots,\vec{v}_n)$ is not transverse for $n\ge 4$.
  \begin{proof}
    Assume $P_+(\vec{v}_1,\dots,\vec{v}_n)$ is transverse. Then, by Lemma \ref{lem:positive_segment_pair} the symplectic products $\w(\vec{v}_i,\vec{v}_j)$ have the same sign for all $i<j-1$. Moreover, if $n>4$ the segment $(1-t)\vec{v}_3 + t \vec{v}_4$ is transverse to the segment $(1-t)\vec{v}_n + t \vec{v}_1$, so $\w(\vec{v}_3,\vec{v}_n)$ has the same sign as $\w(\vec{v}_3,\vec{v}_1) = -\w(\vec{v}_1,\vec{v}_3)$, a contradiction. The case $n=4$ can be normalized to have the matrix of symplectic products $\Omega$ (equation (\ref{eqn:OmegaMatrix})) and the finitely many cases which can be checked by hand are all non-transverse (see Prop. \ref{prop:quadrilaterals}).
  \end{proof}
\end{prop}

Combining Propositions \ref{prop:pplusnottransverse} and \ref{prop:pmcontractible}, we obtain

\begin{cor}
  A transverse Legendrian polygon is non-contractible.
\end{cor}

For non-contractible Legendrian polygons, we now give a simple criterion for transversality.

\begin{prop}\label{prop:pmtransversecondition}
  The Legendrian polygon $P_-(\vec{v}_1,\dots,\vec{v}_n)$ is transverse if and only if $\w(\vec{v}_i,\vec{v}_j)>0$ for all $i < j - 1$, or $\w(\vec{v}_i,\vec{v}_j) < 0$ for all $i < j - 1$.
  \begin{proof}
    We start with the reverse direction. The transversality for segments $(1-t)\vec{v}_i + t \vec{v}_{i+1}$ and $(1-t)\vec{v}_j + t \vec{v}_{j+1}$ with $i < j - 1$ follows at once from Lemma \ref{lem:positive_segment_pair}. It remains to show that segments $(1-t)\vec{v}_i + t \vec{v}_{i+1}$ are transverse to the final segment $(1-t)\vec{v}_n - t \vec{v}_1$. This follows from another application of Lemma \ref{lem:positive_segment_pair} since $\w(\vec{v}_i,-\vec{v}_1) = \w(\vec{v}_1,\vec{v}_i)$ and $\w(\vec{v}_{i+1},-\vec{v}_1) = \w(\vec{v}_1,\vec{v}_{i+1})$.
    
    For the forwards direction, if the polygon is transverse then the non-adjacent vertices are non-incident and so $\w(\vec{v}_i,\vec{v}_j) \neq 0$ for $i < j-1$. Using Lemma \ref{lem:positive_segment_pair} we find that $\w(\vec{v}_i,\vec{v}_j)$ all have the same sign.
  \end{proof}
\end{prop}

The common sign of $\w(\vec{v}_i,\vec{v}_j)$ for a transverse Legendrian polygon $P_-(\vec{v}_1,\dots,\vec{v}_n)$ does not depend on the choice of representatives for the vertices $\vec{v}_1,\dots,\vec{v}_n$. Indeed, multiplication by a positive scalar does not change this sign, and simultaneously changing the sign of all $\vec{v}_j$ also does not change this sign. We will call the corresponding classes of Legendrian polygons \emph{positive-transverse} and \emph{negative-transverse} according to this sign. Reversing the cyclic order exchanges positive-transverse and negative-transverse polygons.
  
  We now recall the definition of the space of generic Lagrangian configurations from \cite{ConleyOvsienko} and relate it to the space of transverse Legendrian polygons.
  
  \begin{defn}
    The moduli space of generic Lagrangian configurations $\mathcal{L}_{2,n}$ is the space of $n$-tuples $(p_1,\dots,p_n)$ in $\Proj(\V)$ such that a pair $p_i,p_j$ is incident if and only if it is cyclically adjacent, up to the action of $\PSp(\V)$.
  \end{defn}
  
  The projection sending a transverse Legendrian $n$-gon to its tuple of vertices is two-to-one, each fiber consisting of a polygon with positive symplectic products of vertices and one with negative products. It is surjective onto $\mathcal{L}_{2,n}$ since is it always possible to choose representative vectors $(\vec{v}_1,\dots,\vec{v}_n)$ for the points $(p_1,\dots,p_n)\in \mathcal{L}_{2,n}$ such that $\omega(\vec{v}_i,\vec{v}_j) > 0$ for all $1 < i+1<j\le n$. Then, the resulting polygon $P_-(\vec{v}_1,\dots,\vec{v}_n)$ is transverse and projects to $(p_1,\dots,p_n)$.
  
  As a consequence of Proposition 2.7 of \cite{ConleyOvsienko}, these two moduli spaces are smooth manifolds of dimension $2(n-5)$ when $n\geq 5$.
  
  If we consider instead generic Legendrian $n$-gons, then the map $\mathscr{P}_n \rightarrow \mathcal{L}_{2,n}$ is no longer surjective, but it is a $2^n$-sheeted covering map of its image which is an open submanifold. We conclude that $\mathscr{P}_n$ is a manifold of dimension $2(n-5)$.

    \begin{prop}\label{prop:quadrilaterals}
      The moduli space of transverse Legendrian $4$-gons is a pair of points.
      \begin{proof}
          Let $P$ be a Legendrian $4$-gon. Then, up to rescaling, its vertices form a basis in which the symplectic form is the standard
  \[\begin{pmatrix}
  0 & 0 & 1 & 0\\
  0 & 0 & 0 & 1\\
  -1 & 0 & 0 & 0\\
  0 & -1 & 0 & 0\end{pmatrix}.\]
  In this basis, the pointwise stabilizer of the four projective vertices in $\PSp(V)$ is the collection of diagonal matrices. Naming the four basis vectors $\vec{e}_1,\vec{e}_2,\vec{e}_3,\vec{e}_4$, the eight possible Legendrian quadrilaterals are $P_\pm(\vec{e}_1,\vec{e}_2,\vec{e}_3,\vec{e}_4)$, $P_\pm(\vec{e}_1,\vec{e}_2,-\vec{e}_3,\vec{e}_4)$, $P_\pm(\vec{e}_1,\vec{e}_2,\vec{e}_3,-\vec{e}_4)$, $P_\pm(\vec{e}_1,\vec{e}_2,-\vec{e}_3,-\vec{e}_4)$, where we normalized the first two basis vectors by applying diagonal matrices with diagonal entries $(-1,1,1,-1)$ and $(1,-1,-1,1)$. By Proposition \ref{prop:pplusnottransverse} and Lemma \ref{lem:positive_segment_pair}, we conclude that out of these only $P_-(\vec{e}_1,\vec{e}_2,\vec{e}_3,\vec{e}_4)$ and $P_-(\vec{e}_1,\vec{e}_2,-\vec{e}_3,-\vec{e}_4)$ are transverse.
      \end{proof}
    \end{prop}

  We now give a definition of the \emph{Maslov index}, a classical invariant of triples of real Lagrangians taking values in the integers.
  Let $L_1,L_2$ be transverse Lagrangians. Define an anti-symplectic involution $\sigma_{L_1,L_2} := I \oplus -I$ according to the splitting $\V = L_1 \oplus L_2$. 
  
  \begin{defn}
    Let $L_1,L_2,L_3$ be pairwise transverse Lagrangians in $\V$. The \emph{Maslov form} is the symmetric nondegenerate bilinear form on $L_2$ defined by
    \[b_{L_1,L_2,L_3}(\vec{u},\vec{v}) = \omega(\vec{u}, \sigma_{L_1,L_3} \vec{v})|_{L_2}.\]
    
    The \emph{Maslov index} $\M(L_1,L_2,L_3)$ is the signature of $b_{L_1,L_2,L_3}$, that is, $\M(L_1,L_2,L_3) = \frac{k_+ - k_-}{2}$ where $k_+$ is the number of positive eigenvalues and $k_-$ is the number of negative eigenvalues of $b_{L_1,L_2,L_3}$.
  \end{defn}

  The Maslov index takes its values in the set $\{-1,0,1\}$ and is a complete $\PSp(2n,\bR)$ invariant of triples of pairwise-transverse Lagrangians.

\section{Legendrian polygons and piecewise circular curves}\label{sec:circles}
In this section, we establish a dictionary between the contact projective geometry of $\Proj(\V)$ and the geometry of co-oriented circles in the $2$-sphere.

Let $(\V,\omega)$ be a $4$-dimensional real symplectic vector space as before. Let $J$ be a complex structure on $\V$ which is anti-symplectic, meaning $J$ is a linear automorphism of $\V$ such that $J^2 = -I$ and $J^*\omega = -\omega$.

As the name suggests, $J$ endows $\V$ with the structure of a complex vector space where multiplication by $i$ is given by the automorphism $J$. Moreover, since $J$ is anti-symplectic, the bilinear form
\[\omega_J(\vec{u},\vec{v}):= \omega(\vec{u},J\vec{v})\]
 is also symplectic.
Together, $\omega$ and $\omega_J$ define a $\bC$-linear symplectic form $\w_\bC$ on $\V$:
\[\omega_\bC := \omega - i \omega_J.\]

Consider the projection
\[\pi : \Proj_\bR(\V) \rightarrow \Proj_\bC(\V)\]
which maps the real projectivisation of a vector $\vec{v}\in \V$, denoted by $[\vec{v}]_\bR$, to its complex projectivisation $[\vec{v}]_\bC$. This projection $\pi$ is a circle bundle, and it is related to the Hopf fibration $p : \mathbb{S}^3 \rightarrow \mathbb{S}^2\cong \mathbb{CP}^1$ by $p = \pi \circ \iota$, where $\iota$ is the covering map $\mathbb{S}^3 \rightarrow \RP^3$.

\begin{defn}
  The \emph{space of oriented contact elements} of $\Proj_\bC(\V)$ is the spherical cotangent bundle $\mathbb{S}(T^*\Proj_\bC(\V))$. It is naturally a contact manifold with the contact distribution given by the tautological $1$-form (or Liouville $1$-form).
\end{defn}

The tangent space at a point $[\vec{v}]_\bC \in \Proj_\bC(\V)$ identifies with the space of $\bC$-linear maps $\Hom_\bC([\vec{v}]_\bC,\V/[\vec{v}]_\bC)$. Similarly, the cotangent space at a point $[\vec{v}]_\bC\in\Proj_\bC(\V)$ identifies with the space of complex linear maps $\Hom_\bC(\V/[\vec{v}]_\bC,[\vec{v}]_\bC)$.

We can evaluate a covector $\alpha \in T^*_{[\vec{v}]_\bC} \Proj_\bC(\V)$ on a vector $X\in T_{[\vec{v}]_\bC}\Proj_\bC(\V)$ by composition:
\[\alpha \circ X \in \Hom([\vec{v}]_\bC,[\vec{v}]_\bC) \cong \bC.\]
Since we are mostly concerned with real manifolds, when $U$ is a complex vector space we use the isomorphism $\alpha \mapsto \Re(\alpha)$ between $\Hom_\bC(U,\bC)$ and $\Hom_\bR(U,\bR)$ to identify complex covectors with real covectors.

In order to identify $\mathbb{S}\Hom_\bC(\V/[\vec{v}]_\bC,[\vec{v}]_\bC)$ with the space of oriented contact elements at $[\vec{v}]_\bC$, we note that an element $\alpha \in \Hom_\bC(\V/[\vec{v}]_\bC,[\vec{v}]_\bC)$ defines a halfspace of the tangent space at $[\vec{v}]_\bC$ by $\Re(\alpha(X))>0$, and this is invariant under multiplication by a positive real.  In figures we depict an oriented contact element as an arrow perpendicular to the boundary of this halfspace and pointing to its interior (Figures \ref{fig:threegraphs}, \ref{fig:transverseoctagon}, and \ref{fig:convex20}).

The following proposition will allow us to interpret the points of $\Proj_\bR(\V)$ as oriented contact elements of $\Proj_\bC(\V)$.

\begin{prop}\label{prop:RP3ContactElements}
  The map
  \begin{align*}
      F : \Proj_\bR(\V) &\rightarrow \mathbb{S}(T^*\Proj_\bC(\V))\\
  [\vec{v}]_\bR &\mapsto (\vec{u}\mapsto \omega_\bC(\vec{u},\vec{v}) \vec{v}),
  \end{align*}
  where the right hand side is viewed as an element of $\Hom_\bC(\V/[\vec{v}]_\bC,[\vec{v}]_\bC)$ up to multiplication by a positive real, is an isomorphism of circle bundles over $\Proj_\bC(\V)$.
  \begin{proof}
    We first verify that the formula $(\vec{u}\mapsto \omega_\bC(\vec{u},\vec{v}))$ defines an element of $\Hom_\bC(\V/[\vec{v}]_\bC,[\vec{v}]_\bC)$.
    
    Since $\omega_\bC$ is $\bC$-bilinear, the map $\vec{u}\rightarrow \omega_\bC(\vec{u},\vec{v})\vec{v}$ is $\bC$-linear. Moreover, it is well-defined on $\V/[\vec{v}]_\bC$ because
    $\omega_\bC(\vec{u}+\vec{v},\vec{v}) = \omega_\bC(\vec{u},\vec{v})$
    and
    $\omega_\bC(\vec{u}+J\vec{v},\vec{v}) = \omega_\bC(\vec{u},\vec{v})$.
    
    For any $t\in \bR^*$, we have
    \[\omega_\bC(\vec{u},t\vec{v})t\vec{v} = t^2\omega_\bC(\vec{u},\vec{v})\vec{v},\]
    so the map is well-defined on the quotients.
    
    Identifying $\V\cong \bC^2$, we can locally trivialize the bundle $\Proj_\bR(\V)$ over the affine patch  $\left\{\begin{bmatrix}z\\1\end{bmatrix}\in\CP^1\right\} \subset \Proj_\bC(\V)$ by
    
    \begin{align*}
        \varphi : \bC \times S^1 &\rightarrow \Proj_\bR(\bC^2)\\
        (z,e^{i\theta}) &\mapsto \begin{bmatrix}z e^{\frac{i\theta}{2}}\\e^{\frac{i\theta}{2}}\end{bmatrix}.
    \end{align*}
    Similarly, we can locally trivialize $\Hom_\bC(\V/[\vec{v}]_\bC,[\vec{v}]_\bC)$ by
    \begin{align*}
        \psi : \bC \times S^1 &\rightarrow \mathbb{S}\left(\Hom_\bC(\V/[\vec{v}]_\bC,[\vec{v}]_\bC)\right)\\
        (z,e^{i\theta}) &\mapsto \left[\vec{u}\mapsto e^{i\theta}\omega_\bC\left(\vec{u},\begin{pmatrix}z\\1\end{pmatrix}\right)\begin{pmatrix}z\\1\end{pmatrix}\right].
    \end{align*}
    In this pair of trivializations, the map $F$ is the identity and hence is continuous and an isomorphism on each fiber. Similar trivilizations on the affine patch $\left\{\begin{bmatrix}1\\z\end{bmatrix}\in\CP^1\right\}$ yield the same result, and we conclude that $F$ is a bundle isomorphism.
  \end{proof}
\end{prop}

The group
\[\PSp(2,\bC) := \{M\in \PSp(\V) ~|~ MJ=JM\}\]
acts transitively on the base $\Proj_\bC(\V)$ and also acts on the total spaces $\Proj_\bR(\V)$ and $\mathbb{S}(T^*\Proj_\bC(\V))$ of the two circle bundles. The action on $\mathbb{S}(T^*\Proj_\bC(\V))$ when viewed as a $\Hom$-bundle is by conjugation, that is, if $f\in \mathbb{S}\Hom_\bC(\V/[\vec{v}]_\bC, [\vec{v}]_\bC)$, then
\[M \cdot f = MfM^{-1} \in \mathbb{S}\Hom_\bC(\V/[M\vec{v}]_\bC, [M\vec{v}]_\bC)\]
is the image of $f$ by the action of $A\in\PSp(2,\bC)$. The isomorphism in Proposition \ref{prop:RP3ContactElements} intertwines the two actions on the total spaces.

\begin{prop}
  The map $F : \Proj_\bR(\V) \rightarrow \mathbb{S}(T^*\Proj_\bC(\V))$ is an isomorphism of contact manifolds.
  \begin{proof}
    We will use the same trivializations $\varphi,\psi$ as in the proof of Proposition \ref{prop:RP3ContactElements}.
    
    In the coordinate chart $\varphi$, a path $\gamma(t) = (z(t),e^{i\theta(t)})$ is tangent to the contact distribution of $\Proj_\bR{\bC^2}$ if and only if
    \[\w_\bR(\gamma(t),\gamma'(t)) = \omega_\bR\left(\begin{pmatrix}z(t) e^{\frac{i\theta(t)}{2}}\\e^{\frac{i\theta(t)}{2}}\end{pmatrix},
    \begin{pmatrix}
    z'(t)e^{\frac{i\theta(t)}{2}} + z(t)\frac{i\theta'(t)}{2}e^{\frac{i\theta(t)}{2}}\\
    \frac{i\theta'(t)}{2}e^{\frac{i\theta(t)}{2}}
    \end{pmatrix}\right) = 0,\]
    which simplifies to
    \[\Re\left(e^{i\theta(t)}\omega_\bC\left(\begin{pmatrix}z(t)\\1\end{pmatrix},
    \begin{pmatrix}
    z'(t)\\ 0 \end{pmatrix}\right)\right) = 0.\]
    
    In the coordinate chart $\psi$ for $\mathbb{S}(T^*\CP^1)$, the path $\gamma(t) = (z(t),e^{i\theta(t)})$ is in the contact distribution of if and only if
    \[\psi(\gamma(t))\left(d\pi\left(\left.\frac{\dif}{\dif s}\right|_{s=t}\psi(\gamma(s))\right)\right) = 
    \Re\left(e^{i\theta(t)}\omega_\bC\left(\begin{pmatrix} z'(t)\\0\end{pmatrix},\begin{pmatrix}z(t)\\1\end{pmatrix}\right)\right) = 0.\]
    
    Since these equations define the same plane, and in these two charts the isomorphism $F$ is the identity, we conclude that $F$ is a contactomorphism.
  \end{proof}
\end{prop}

Now that we have established the contact isomorphism between $\Proj(\V)$ and the space of oriented contact elements to the sphere, we can begin to explore the geometric interpretations. The object which corresponds to a Legendrian line (of Lagrangian subspace) is a co-oriented circle.

\begin{defn}
  A \emph{co-oriented circle} in $\Proj_\bC(\V)$ is a round circle $C$ of $\Proj_\bC(\V)$, possibly degenerated to a point, together with a choice of disk bounded by $C$ when nondegenerate.
\end{defn}

We say that two co-oriented circles are tangent if the circles are tangent and the orientations match at the tangency point. If one of the circles is degenerate, they are considered tangent if the point circle is contained in the other circle.

\begin{lem}
  Let $L\subset \V$ be a Lagrangian subspace for $\omega$. If $L$ is not a complex subspace, then $\pi(L)$ is a circle. The two connected components of $\Proj_\bC(\V) \setminus \pi(L)$ are distinguished by the Maslov index :
  \[\{\ell \in \Proj_\bC(\V) ~|~ \M(L,\ell,JL) = 1\}\]
  and
  \[\{\ell \in \Proj_\bC(\V) ~|~ \M(L,\ell,JL) = -1\}.\]
  \begin{proof}
    Since $L$ is not a complex subspace, $L\oplus JL = \V$ defines a splitting of the vector space $\V$ into a pair of Lagrangians.
    
    The anti-symplectic, anti-linear involution $\sigma_{L,JL} = I \oplus -I$ defines a \emph{real structure} on the complex vector space $\V$. When endowed with this real structure, the complex projectivization $\pi(L)$ identifies with the real projective line $\Proj_\bR(\V)\subset \Proj_\bC(\V)$, which is a (generalized) circle.
    
    The involution $\sigma_{L,JL}$ exchanges the two connected components of the complement $\Proj_\bC(\V) \setminus \pi(L)$. The Maslov index $\M(L,\ell,JL)$ is the signature of the symmetric bilinear form $b_{L,\ell,JL} = \omega(\cdot,\sigma_{L,JL}(\cdot))$ restricted to $\ell$.
    
    Note that $J$ is an automorphism of the bilinear form $b_{L,\ell,JL}$:
    \[\omega(J\vec{u},\sigma_{L,JL}J\vec{v}) = \omega(J\vec{u},-J\sigma_{L,JL}\vec{v}) = \omega(\vec{u},\sigma_{L,JL}\vec{v}).\]
    
    Moreover, $\vec{v}$ is always orthogonal to $J\vec{v}$ for the form
    $b_L,J_L$. Hence, if $\ell$ is the complex line spanned the vector $\vec{v}$, $b_{L,\ell,JL}$ is either definite or identically zero. If it is identically zero, then $\omega(\vec{v},\sigma_{L,JL}\vec{v}) = 0$ and $\omega(J\vec{v},\sigma_{L,JL}\vec{v}) = 0$ and the fact that $\ell = \Span_\bR(\vec{v},J\vec{v})$ is a Lagrangian imply that $\sigma_{L,JL}(\vec{v})\in \ell$, and hence $\ell \in \pi(L)$.
    
    By continuity of the bilinear form $b_{L,\ell,JL}$ we conclude that the two complementary components of $\pi(L)$ consist of where it is positive definite and negative definite.
   \end{proof}
\end{lem}

Let $L\subset \V$ be a Lagrangian subspace for the symplectic form $\omega$. If $L\neq JL$, let $C_L$ denote the co-oriented circle $\pi(L)$ with the co-orientation given by the disk of complex lines $U\subset \V$ such that $\M(L,U,JL) = -1$. If $L=JL$, then $L$ is a complex subspace and we let $C_L = \pi(L)$ which we interpret as a zero radius circle. Then, $L\leftrightarrow C_L$ is a bijection between the set of Lagrangians in $\V$ and the set of co-oriented circles in $\Proj_\bC(\V)$.

\subsection{In coordinates}\label{sec:coordinates}
Concretely, in a basis $\vec{e}_1,\vec{e}_2,\vec{e}_3,\vec{e}_4$ for which the symplectic form is given by \[\Omega = \begin{pmatrix} 0 & 0 & 0 & 1\\0 & 0 & -1 & 0\\0 & 1 & 0 & 0\\ -1 & 0 & 0 & 0\end{pmatrix},\]
a dense open set of Lagrangians is parametrized by the affine patch
\[L(a,b,c) = \begin{pmatrix} 1 & 0 \\0 & 1\\a & b+c\\ b-c & -a\end{pmatrix}.\]

Choosing the compatible complex structure with matrix
\[ J = \begin{pmatrix} 0 & 1 & 0 & 0\\-1 & 0 & 0 & 0\\0 & 0 & 0 & -1\\ 0 & 0 & 1 & 0\end{pmatrix},\]
the complex structure acts on Lagrangians $L(a,b,c)$ by reversing the sign of $c$, and so the Lagrangians in the affine patch which are complex lines are precisely $L(x,y,0)$. The pair $(\vec{e}_1,\vec{e}_3)$ forms a complex basis for $\V$ in which the complex symplectic form $\omega_\bC$ has matrix $\begin{pmatrix}0 & -i\\i & 0\end{pmatrix}$. In this basis, we have
\[L(x,y,0) = \begin{pmatrix} 1 & 0 \\0 & 1\\x & y\\ y & -x\end{pmatrix} = \begin{pmatrix}1 \\ x + iy\end{pmatrix}.\]

The Maslov index $\M(L(a,b,c),L(x,y,0),JL(a,b,c))$ is the index of the bilinear form
\[\begin{pmatrix}\frac{(x-a)^2 + (y-b)^2-c^2}{c} & 0\\ 0 & \frac{(x-a)^2 + (y-b)^2-c^2}{c}\end{pmatrix},\]
which is negative definite if and only if either: $c<0$ and $(x,y)$ is inside of the circle of radius $|c|$ centered at $a+bi$, or $c>0$ and $(x,y)$ is outside of that circle. Thus, $L(a,b,c)$ represents the circle of radius $|c|$ centered at $a+bi$ co-oriented towards the inside if $c>0$ and towards the outside if $c<0$.

\begin{figure}[h]
    \centering
    \includegraphics[width=.6\textwidth]{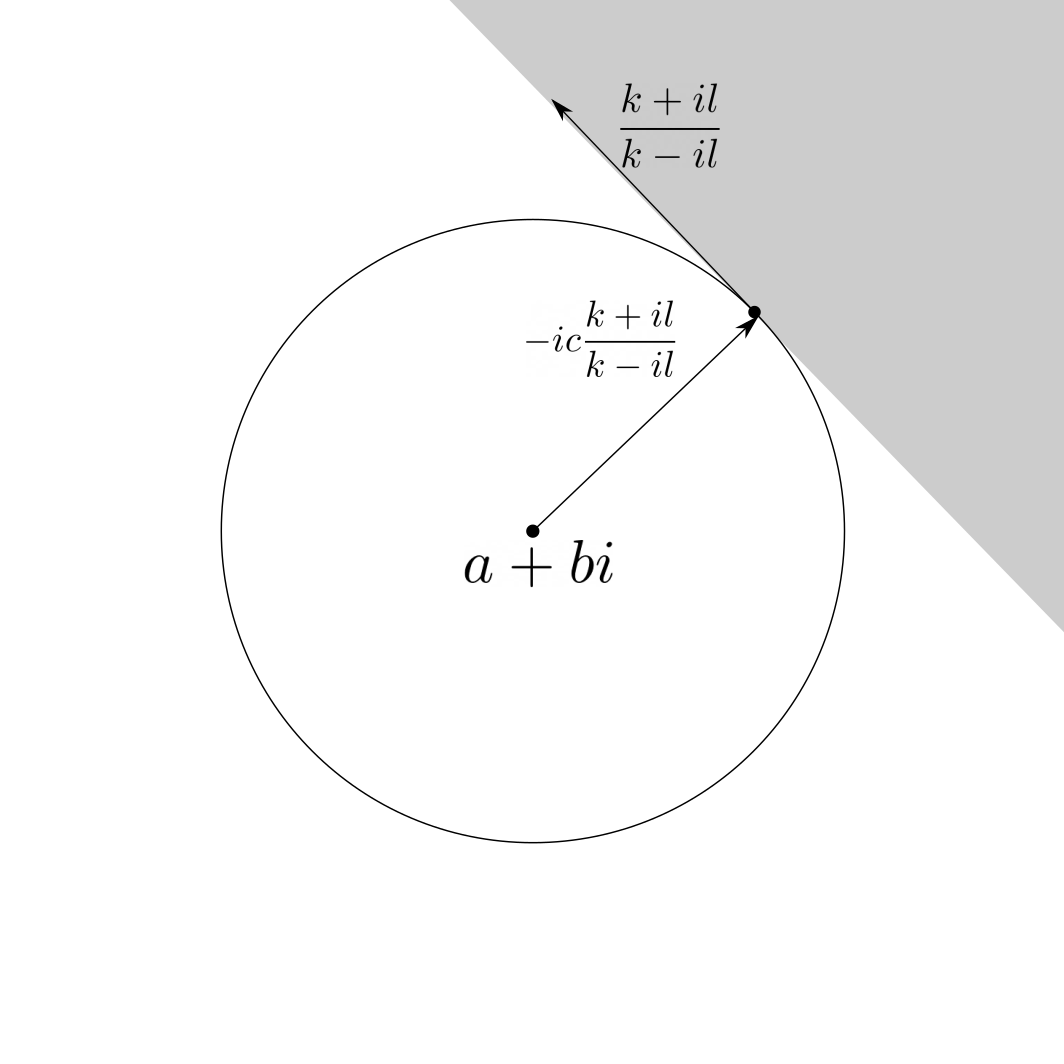}
    \caption{An oriented contact element to the circle corresponding to $L(a,b,c)$ with $c>0$, in local coordinates.}
    \label{fig:covec_circle}
\end{figure}

\begin{lem}
  Let $\vec{v}\in \V$ and $L\subset \V$ be a Lagrangian. Then, $\vec{v}\in L$ if and only if $\pi(\vec{v})\in C_L$ and the oriented contact element to $C_L$ at $\pi(\vec{v})$ is given by $[\vec{v}]_\bR$ through the isomorphism $F$ of \hyperref[prop:RP3ContactElements]{Proposition \ref{prop:RP3ContactElements}}.
  \begin{proof}
    Using homogeneity under the $\PSp(2,\bC)$ action, it would be enough to show this for a single $\vec{v}\in\V$. Since the calculation is simple enough, we show it for any point which projects to a standard affine patch in $\CP^1$.
    
    In affine coordinates as above, we consider
    \[[\vec{v}]_\bR = \begin{bmatrix}k\\l\\k a + l (b+c)\\k(b-c)-l a\end{bmatrix}_\bR \subset L(a,b,c).\]
    
    This point in real projective space projects to
    \[\pi([\vec{v}]_\bR) = [\vec{v}]_\bC = \begin{bmatrix} k - i l\\k(a+(b-c)i)+l(b+c) - ai\end{bmatrix} = \begin{bmatrix} 1\\a + b i - i c\frac{k+il}{k-il}\end{bmatrix}\in \CP^1.\]
    
    This is a parametrization by $[k,l]\in\RP^1$ of the circle corresponding to $L(a,b,c)$ where $[\vec{v}]_\bC$ is the unique point on the circle of radius $|c|$ centered at $a+bi$ which is in the direction of $-i\frac{(k+il)^2}{|k+il|^2}$ if $c>0$, and in the opposite direction if $c<0$.
    
    The tangent vector at $[\vec{v}]_\bC$ corresponding to $z \in \bC$, as a linear map, is represented by
    \[X_z = \vec{v}_\bC^* \otimes \begin{pmatrix} 0\\(k-il)z\end{pmatrix} \in \Hom_\bC([\vec{v}]_\bC , \bC^2/[\vec{v}]).\]
    
    Evaluating the covector $F([\vec{v}]_\bR)$ at $X_z$ yields
    \[\omega_\bC\left(\begin{pmatrix}0\\(k-il)z\end{pmatrix},\vec{v}\right) = i(k-il)^2z.\]
    
    We conclude that the halfspace of the tangent space at $[\vec{v}]_\bC$ defined by the covector $F([\vec{v}]_\bR)$ is
    \[\{z \in \bC ~|~ \Re(zi(k-il)^2) > 0\} = \{z \in \bC ~|~ \Re(z\overline{(-i(k+il)^2)}) > 0\}.\]
    This halfspace is bounded by the line spanned by $(k + il)^2 = |k+il|^2\frac{k+il}{k-il}$ and contains $-i(k+il)^2$.
    
    Therefore, if $c>0$, the circle corresponding to $L(a,b,c)$ is oriented outwards and the halfspace of the tangent space at $a+bi-ic\frac{k+il}{k-il}$ determined by $F([\vec{v}]_\bR)$ is also oriented outwards. If $c<0$, the circle corresponding to $L(a,b,c)$ is oriented inwards and the halfspace of the tangent space is also oriented inwards (see Fig. \ref{fig:covec_circle}).
    
    We conclude that the oriented contact elements represented by lines in $\bR^4$ contained in the Lagrangian $L(a,b,c)$ are precisely the oriented contact elements of the co-oriented circle corresponding to $L(a,b,c)$, proving the proposition.
  \end{proof}
\end{lem}

\begin{cor}
  Let $[\vec{u}]_\bR, [\vec{v}]_\bR\in \Proj_\bR(\V)$. Then, $\w(\vec{u}, \vec{v}) = 0$ if and only if $F([\vec{u}]), F([\vec{v}])$ are oriented contact elements tangent to a common co-oriented circle.
\end{cor}

In the coordinates we are using, the embedding $\PSL(2,\bC)\cong \PSp(2,\bC) \hookrightarrow \PSp(4,\bR)$ induced by forgetting the imaginary part of the complex symplectic form $\omega_\bC$ is
\[\begin{pmatrix}
a & b\\
c & d
\end{pmatrix}
\mapsto
\begin{pmatrix}
 \Re(a) & \Im(a) & \Re (b) & -\Im (b) \\
 -\Im (a) & \Re (a) & -\Im (b) & -\Re (b) \\
 \Re (c) & \Im (c) & \Re (d) & -\Im (d) \\
 \Im (c) & -\Re (c) & \Im (d) & \Re (d) \\
\end{pmatrix}.\]

These matrices act by M\"obius transformations on the unit tangent bundle of $S^2$.

Let $r\in \bR$ and define
  \[T_r := \begin{pmatrix}
  1 & 0 & 0 & 0\\
  0 & 1 & 0 & 0\\
  0 & r & 1 & 0\\
  -r & 0 & 0 & 1
  \end{pmatrix}.\]
  Then, $T_r(L(a,b,c)) = L(a,b, c + r)$. In geometric terms, the transformation $T_r\in \PSp(4,\bR)$ maps any oriented circle in the affine patch $L$ to the circle with the same center and radius increased by $r$. Its action on piecewise circular curves is therefore the radial translation of parameter $r$.

\begin{prop}
  M\"obius transformations and radial translations $T_r$ generate the Lie group $\PSp(4,\bR)$.
  \begin{proof}
     The Proposition follows by a Lie algebra computation using the matrix expressions of the two subgroups and the fact that $\PSp(4,\bR)$ is connected.
  \end{proof}
\end{prop}

Therefore, the classification problem of piecewise circular $n$-gons up to M\"obius transformations and radial translations reduces to the classification of Legendrian $n$-gons in $\RP^3$ modulo the action of $\PSp(4,\bR)$.

We can now justify the interpretation of the Maslov index as a generalization of nestedness of circles mentioned in the introduction. The group $\PSp(4,\bR)$ acts transitively on pairs of transverse Lagrangians, and transitively on pairwise-transverse triples of Lagrangians which have the same Maslov index. Therefore, any triple with Maslov index $+1$ is equivalent to the triple $L(0,0,0)$, $L(0,0,1)$, $L(0,0,2)$ corresponding to three nested circles.

\begin{figure}[ht]
    \centering
    \begin{subfigure}[b]{0.45\textwidth}
         \centering
         \includegraphics[width=\textwidth]{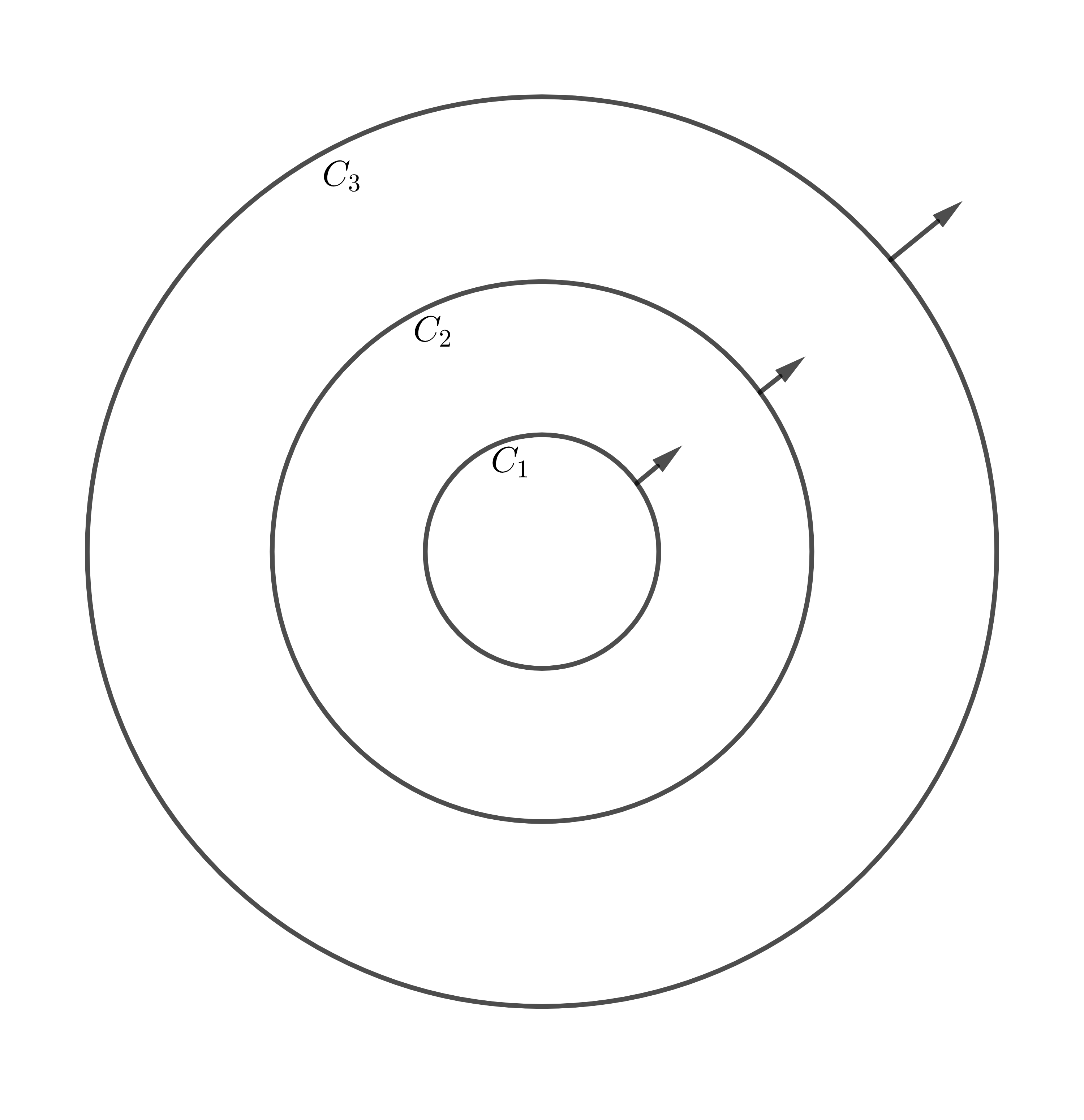}
         \caption{$\M=1$}
     \end{subfigure}
     \begin{subfigure}[b]{0.45\textwidth}
         \centering
         \includegraphics[width=\textwidth]{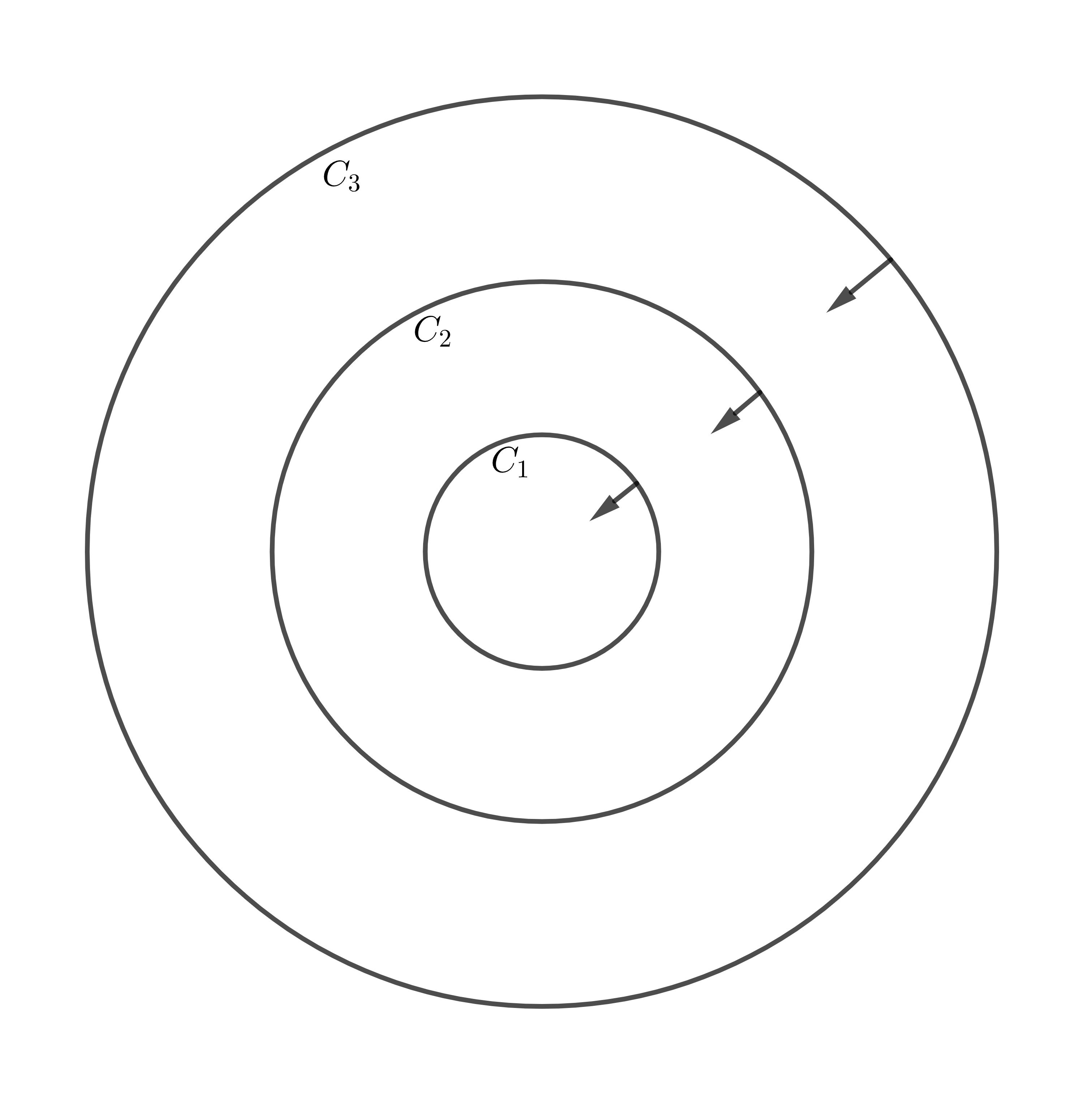}
         \caption{$\M=-1$}
     \end{subfigure}\\
     \begin{subfigure}[b]{0.9\textwidth}
         \centering
         \includegraphics[width=\textwidth]{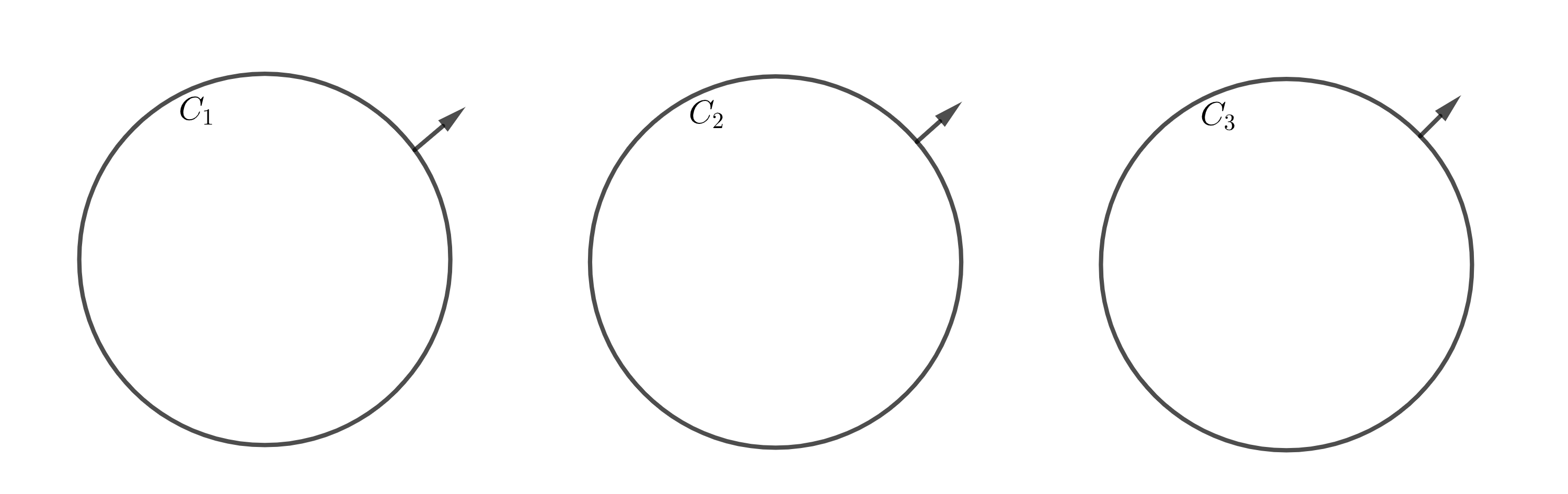}
         \caption{$\M=0$}
     \end{subfigure}
        \caption{The Maslov index of a triple of co-oriented circles.}
        \label{fig:threecircles}
\end{figure}

In other words, any configuration of three co-oriented circles which are pairwise non tangent can be brought by a sequence of M\"obius transformations and radial translations to exactly one of the configurations depicted in Figure \ref{fig:threecircles}, and Maslov index $1$ corresponds to nested circles with increasing radius (or equivalently, decreasing curvature).

\begin{defn}
  We say that a transverse Legendrian polygon and its corresponding piecewise circular curve have \emph{decreasing curvature} if every triple of non-adjacent segments is a triple of Lagrangians with Maslov index $1$.
\end{defn}


Recall that we defined a \emph{simple} piecewise circular curve by the property that every radial translate is simple.

\begin{prop}\label{prop:simpleifftransverse}
  A piecewise circular curve is simple if and only if the corresponding Legendrian polygon is transverse.
  \begin{proof}
    A self-intersection of a piecewise circular curve can be interpreted as a zero-radius circle which is tangent to two oriented contact elements of the curve lying on non-adjacent pieces. This means that there exists a Lagrangian containing two points in $\Proj(\V)$ lying on non-adjacent edges of the corresponding Legendrian polygon. In other words, those two points are incident, so the Legendrian polygon is not transverse.
    
    Conversely, if a Legendrian polygon is non-transverse, then there exist two points on non-adjacent edges which are incident. The Lagrangian spanned by these two points corresponds to a circle tangent to two non-adjacent circular pieces of the piecewise circular curve. Applying a M\"obius transformation if needed, we may assume that this circle is not a line. Then, applying a radial translation by the opposite of its radius brings that circle to a zero radius circle, producing a translate with a self-intersection.
  \end{proof}
\end{prop}

\subsection{The dictionary}
From the above analysis, we deduce the dictionary in Table \ref{tab:dict} between the geometry of the symplectic vector space $(\V,\omega)$ and the geometry of circles in the Riemann sphere $\Proj_\bC(\V)$.

\begin{table}
{\renewcommand{\arraystretch}{1.6}%
\begin{center}
    \begin{tabular}{|c|c|}
    \hline
    $(\V,\omega)$ & $\Proj_\bC(\V)$\\
    \hline
    \hline
    Point in $\Proj_\bR(\V)$     &  Oriented contact element\\
    \hline
    Lagrangian subspace of $\V$     &  Co-oriented circle\\
    \hline
    Isotropic flag in $\V$ & Pointed co-oriented circle\\
    \hline
    Non-transverse Lagrangians & Tangent co-oriented circles\\
    \hline
    Incident points in $\Proj_\bR(\V)$ & Oriented contact elements\\
     & to a common co-oriented circle\\
    \hline
    Legendrian polygon & Co-oriented\\ & piecewise circular curve\\
    \hline
    Moduli space of generic  & Moduli space of generic \\
    Legendrian $n$-gons & circular $n$-gons\\
    \hline
    Transverse Legendrian $n$-gon & Simple circular $n$-gon\\
    \hline
    \end{tabular}
\end{center}}
\caption{Dictionary between Legendrian polygons and piecewise circular curves}
\label{tab:dict}
\end{table}

We can translate the results of the previous sections with this dictionary as follows.

\begin{cor}
  The moduli space of simple circular $n$-gons and the moduli space of generic circular $n$-gons are a smooth manifolds of dimension $2(n-5)$ when $n\geq 5$.
\end{cor}

\begin{cor}
  The moduli space of simple circular quadrilaterals is a pair of points (Proposition \ref{prop:quadrilaterals}).
\end{cor}


\section{Positivity} \label{sec:Positivity}
  \subsection{Flag positivity}
  In this section, we describe positivity in the oriented flag manifold $G/B_0$ of the group $G = \Sp(\V,\omega) \cong \Sp(4,\bR)$. For a more general perspective on positivity in oriented flag manifolds see ~\cite{burelle2018schottky}.
  
  Let $F$ be a flag in $\V$. We denote by $F^{(k)}$ the $k$-dimensional part of $F$.
  
  Note that the symplectic vector space $\V$ is canonically oriented by the volume form $\omega\wedge\omega$. Whenever $U,W \subset \V$ are oriented vector subspaces of $\V$, the equality $U=W$ will mean that the subspaces are equal and the orientations agree.
  
  \begin{defn}
    An \emph{isotropic flag} in $\V$ is a full flag $F$ such that $F^{(4-k)}=F^{(k)\perp}$ for all $k$.
  \end{defn}
  
  Let $F$ be an isotropic flag and choose orientations on each $F^{(k)}$. The orientations on $F^{(1)}$ and $F^{(3)}$ induce an orientation on $F^{(3)}/F^{(1)}$ in the following way: if $(\vec{e}_1,\vec{e}_2,\vec{e}_3)$ is an oriented basis of $F^{(3)}$ such that $\vec{e}_1$ is an oriented basis of $F^{(1)}$, then $(\vec{e}_2,\vec{e}_3)$ gives a well-defined orientation on the quotient. We say that the orientations on $F^{(1)}$ and $F^{(3)}$ are \emph{compatible} if this orientation on $F^{(3)}/F^{(1)}$ matches the orientation induced on this quotient by the symplectic form $\omega$. 
  
  \begin{defn}
    An \emph{oriented isotropic flag} in $\V$ is an isotropic flag $F$ together with a choice of orientations satisfying:
    \begin{itemize}
        \item The orientations on $F^{(1)}$ and $F^{(3)}=F^{(1)\perp}$ are compatible;
        \item The orientation on $F^{(4)}$ is the same as that of $\V$.
    \end{itemize}
  \end{defn}
  
  \begin{rmk}\label{rmk:flagdeterminedby12}
  Because of the condition on compatible orientations, an oriented isotropic flag $F$ is uniquely determined by $F^{(1)}$ and $F^{(2)}$.
  \end{rmk}
  
  We call a basis $E = (\vec{e}_1,\vec{e}_2,\vec{e}_3,\vec{e}_4)$ of $\V$ in which the symplectic form is $\Omega$ (defined in Equation (\ref{eqn:OmegaMatrix})) a \emph{symplectic basis}.
  
  A symplectic basis determines an oriented isotropic flag $F_E = \Span(\vec{e}_1) \subset \Span(\vec{e}_1,\vec{e}_2)\subset \Span(\vec{e}_1,\vec{e}_2,\vec{e}_3) \subset \V$. We now fix such a symplectic basis $E$ and we denote an oriented isotropic flag $F$ by a $4\times4$ matrix such that the first $k$ columns form an oriented basis of $F^{(k)}$. To shorten notation we will also sometimes use the $4\times2$ consisting of only the first two columns, which is equivalent by Remark \ref{rmk:flagdeterminedby12}.
  
  \begin{defn}
     A pair of oriented isotropic flags $F,G$ is called \emph{oriented-transverse} if there exists a symplectic basis $E$ such that $F=F_E$ and $G=F_{\widehat{E}}$.
  \end{defn}
  
  The group $G = \Sp(\V,\omega)$ acts transitively on oriented isotropic flags and the stabilizer of the flag $F_0 := F_E$ is the subgroup $B_+^0$ of upper triangular matrices with positive entries on the diagonal (in the basis $E$). The space of oriented isotropic flags then identifies with the homogeneous space $G/B_+^0$.
  
  The group $G$ also acts transitively on pairs of oriented-transverse oriented isotropic flags $F_E,F_{\widehat{E}}$, with stabilizer diagonal matrices with positive entries.
  

  
  The unipotent radical of the Borel subgroup $B_+$ is the subgroup:
  \[U_+ = \left\{
  \left.\begin{pmatrix}
  1 & a & b & c\\
  0 & 1 & d & a d - b\\
  0 & 0 & 1 & a\\
  0 & 0 & 0 & 1\\
  \end{pmatrix} ~\right|~ a,b,c,d\in \bR\right\}.\]
  
  We define the \emph{opposite symplectic basis} $\widehat{E}$ of $E$ to be $(\vec{e}_4,-\vec{e}_3,\vec{e}_2,-\vec{e}_1)$.
  
  The stabilizer of $F_\infty := F_{\widehat{E}}$ is the subgroup $B_-^0$ of lower triangular matrices with positive entries on the diagonal. The Borel subgroup $B_-$ is opposite to $B_+$, and its unipotent radical is
  \[U_- = \left\{
  \left.\begin{pmatrix}
  1 & 0 & 0 & 0\\
  a & 1 & 0 & 0\\
  b & d & 1 & 0\\
  c & a d - b & a & 1\\
  \end{pmatrix} ~\right|~ a,b,c,d\in \bR\right\}.\]
  
  We say that two oriented isotropic flags $(F_1,F_2)$ are \emph{transverse} if they are in the $\Sp(\V,\omega)$-orbit of the pair $(F_0,F_\infty)$, or equivalently if there is a symplectic basis $E$ such that $F_1 = F_E$ and $F_2 = F_{\widehat{E}}$.
  
  The \emph{positive semigroup} $U_+^{>0}$ is the following subsemigroup of $U_+$:
  \[U_+^{>0} = \left\{
  \left.\begin{pmatrix}
  1 & a & b & c\\
  0 & 1 & d & a d - b\\
  0 & 0 & 1 & a\\
  0 & 0 & 0 & 1\\
  \end{pmatrix} ~\right|~ \begin{array}{c} a>0,~  b>0,~  c>0\\ad-b>0,~ -b^2+abd-cd>0\end{array}\right\},\]
  and similarly
  \[U_-^{>0} = \left\{
  \left.\begin{pmatrix}
  1 & 0 & 0 & 0\\
  a & 1 & 0 & 0\\
  b & d & 1 & 0\\
  c & a d - b & a & 1\\
  \end{pmatrix} ~\right|~ \begin{array}{c} a>0,~  b>0,~  c>0\\ad-b>0,~ -b^2+abd-cd>0\end{array}\right\}.\]
  The inequalities defining the positive semigroups are equivalent to the statement that all minors of the matrix which aren't zero by triangularity are positive (\cite[Theorem 2.8]{pinkus}).
  
  A simple computation shows that these semigroups satisfy $\left(U_\pm^{>0}\right)^{-1} = K U_\pm^{>0} K$, where $K$ is the diagonal matrix $K=\diag(1,-1,1,-1)$.
  
  \begin{defn}
    A triple of pairwise transverse oriented isotropic flags $(F_E,F,F_{\widehat{E}} )$ is \emph{positive} if $F = u F_0$ with $u\in U_-^{>0}$.

    A tuple of flags $(F_1,\dots,F_n)$ is \emph{positive} if every ordered sub-triple $F_i,F_j,F_k$ with $1\le i < j < k \le n$ is positive.
  \end{defn}
  
  In a positive triple, all pairs are oriented-transverse. Additionally, positivity satisfies the following:
  \begin{enumerate}
      \item If $(F_1,F_2,F_3)$ is positive, then $(F_2,F_3,-F_1)$ is positive.
      \item If $(F_1,F_2,F_3)$ is positive, then $(F_3,F_2,F_1)$ is not positive.
      \item If $(F_1,F_2,F_3)$ and $(F_1,F_2,F_4)$ are positive, then $(F_1,F_3,F_4)$ is positive.
  \end{enumerate}
  
  The first property implies that if $(F_1,\dots,F_n)$ is positive, then $(-F_1,\dots,-F_n)$ also is. This motivates the introduction of the following moduli space.
  
  \begin{defn}
     The space $\mathscr{F}^{(n)}$ is the space of $\PSp(\V)$-orbits of $n$-tuples of pairwise transverse oriented isotropic flags up to the diagonal action of $-I$.
     The subspace consisting of positive $n$-tuples is denoted by $\mathscr{F}^{(n)}_{>0}$.
  \end{defn}
  
  As a straightforward consequence of the definitions, we have:
  \begin{prop}\label{prop:positivequadruple}
    A quadruple of oriented isotropic flags $F_0,F_1,F_2,F_\infty$ is positive if and only if $F_1=u_1 F_0$ and $F_2=u_1 u_2 F_0$ with $u_1,u_2\in U_-^{>0}$.
    
    A quadruple of oriented isotropic flags $F_0,F_+,F_\infty,F_-$ is positive if and only if
    $F_+ = u_+ F_0$ and $F_- = -u_- F_0$, where $u_+\in U_-^{>0}$ and $u_-\in KU_-^{>0}K$.
  \end{prop}
  
  The following simple lemma will be useful in the proof of Proposition \ref{prop:polygonpositiveimpliesflagpositive}.
  
  \begin{lem}\label{lem:cyclicsetcharacterization}
    Let $(F_1,\dots,F_n)\in \mathscr{F}^{(n)}$.
    If $(F_i,F_{i+1},F_j)$ is positive for every $1\le i+1 < j \le n$, then $(F_1,\dots,F_n)$ is positive.
    \begin{proof}
      Let $1\le i<j<k\le n$. Since $(F_i,F_{i+1},F_k)$ and $(F_{i+1},F_{i+2},F_k)$ are positive by hypothesis, by property $(1)$ above we find that $(F_k,-F_i,-F_{i+1})$ and $(F_k,-F_{i+1},-F_{i+2})$ are positive. By properties $(3)$ and $(1)$, we conclude that $(F_i,F_{i+2},F_k)$ is positive. Repeating this argument we find that $(F_i,F_{i+m},F_k)$ is positive for all $0 < m < k-i$ and so $(F_i,F_j,F_k)$ is positive.
    \end{proof}
  \end{lem}
  
  Positivity was originally defined for unoriented flags. Let us denote by $\Flag$ the space of unoriented isotropic flags. An $n$-tuple of isotropic flags is called positive if it can be written
  \[F_0,u_1 F_0, u_1 u_2 F_0, \dots, (u_1 u_2\dots u_{n-2})F_0,F_\infty,\]
  where $u_i\in U_-^{>0}$. Let $\Flag^{(n)}$ denote the moduli space of $n$-tuples of pairwise transverse flags up to the action of $\PSp(4,\bR)$.
  
  We will use the following special case of a theorem of Fock and Goncharov  to prove our main result.
  \begin{thm}[{\cite[Theorem 1.5]{FG2006}}]\label{thm:fockgoncharov}
    Positivity of an $n$-tuple of flags is invariant under the action of $\PSp(V)$, and
    the subset $\Flag^{(n)}_{>0} \subset \Flag^{(n)}$ consisting of positive $n$-tuples is a connected component homeomorphic to a ball of dimension $4n-10$.
  \end{thm}
  
  For our oriented isotropic flag setting, we have the same result.
  
  \begin{prop}\label{prop:FockGoncharovComponent}
    The subset $\mathscr{F}^{(n)}_{>0} \subset \mathscr{F}^{(n)}$ consisting of positive $n$-tuples of oriented isotropic flags is a connected component homeomorphic to a ball of dimension $4n-10$.
    \begin{proof}
      The projection $\mathscr{F}^{(n)}_{>0} \rightarrow \Flag^{(n)}$ corresponding to removing the orientations is a covering map. By Theorem \ref{thm:fockgoncharov}, the preimage of the positive component, which is contractible, is a disjoint union of homeomorphic copies.
      
      The definition of a positive $n$-tuples of flags gives a prefered choice of lift to oriented flags, namely the lift
       \[F_0,u_1 F_0, u_1 u_2 F_0, \dots (u_1 u_2\dots u_{n-2})F_0,F_\infty,\]
       now considered as oriented flags. Since the elements $u_i\in U_-$ are uniquely determined, this lift is well-defined and is positive as a tuple of oriented isotropic flags.
       
       This is the only positive lift of this tuple, since we can assume that the lift of $F_0$ is $F_0$ by acting with $\PSp(4,\bR)$, and any pair of unoriented flags $F_0,F_i$ has a unique oriented-transverse lift. Since every positive tuple of oriented isotropic flags is pairwise oriented-transverse, the positive lift is unique.
    \end{proof}
  \end{prop}

    The cited Fock-Goncharov result shows that positive $n$-tuples form a connected component in the space of \emph{generic} $n$-tuples, rather than pairwise-transverse $n$-tuples. However, since positivity is characterized by triples, and a positive triple can only degenerate to a triple containing a non-transverse pair (\cite[Prop. 8.14]{lusztig}), the theorem is true as stated.
    
    The proof of Theorem \ref{thm:fockgoncharov} provides explicit coordinates for $\Flag^{(n)}_{>0}$ and hence for $\mathscr{F}^{(n)}_{>0}$.
    
    \begin{example}\label{ex:parametrization}
      The following map is a parametrization of $\mathscr{F}_{>0}^{(3)}$:
      \begin{align*}
          \bR_{>0} \times \bR_{>0} &\longrightarrow \mathscr{F}^{(3)}_{>0}\\
          (x,y) &\longmapsto (F_0, \left(
\begin{array}{cccc}
 1 & 0 \\
 x+y+\frac{1}{y} & 1 \\
 y & 1 \\
 1 & x+\frac{1}{y} \\
\end{array}
\right), F_\infty ).
      \end{align*}
    \end{example}

  The two-dimensional parts of a triple of oriented isotropic flags are Lagrangians, and so given a triple of flags we obtain a triple of Lagrangians. If the triple is positive, their Maslov index is always $+1$:
  \begin{lem}\label{lem:PositiveMaslov}
    Let $F_1,F_2,F_3$ be a positive triple of oriented isotropic flags. Then, the Maslov index $\M\left(F_1^{(2)},F_2^{(2)},F_3^{(2)}\right)$ is $+1$.
    \begin{proof}
      After acting by an element of $\PSp(4,\bR)$, we may assume
      \[F_1 = \begin{pmatrix} 1 & 0 & 0 & 0\\ 0 & 1 & 0 & 0\\0 & 0 & 1 & 0\\0 & 0 & 0 & 1\end{pmatrix}, \quad
      F_2 = \begin{pmatrix} 1 & 0 & 0 & 0\\ a & 1 & 0 & 0\\ b & d  & 1 & 0\\ c & a d - b & a & 1 \end{pmatrix}, \quad 
      F_3 = \begin{pmatrix} 0 & 0 & 0 & -1\\ 0 & 0 & 1 & 0\\ 0 & -1 & 0 & 0\\ 1 & 0 & 0 & 0\end{pmatrix}\]
      with the matrix for $F_2$ in the semigroup $U_-^{>0}$.
      
      The Maslov bilinear form of the three Lagrangians $F_1^{(2)},F_2^{(2)},F_3^{(2)}$ is then given by 
      \[2\begin{pmatrix} a b - c & b\\b & d \end{pmatrix}.\]
      Since $d>0$ and the determinant $(ab-c)d-b^2 = b(ad-b)-c d>0$, this matrix is positive-definite.
    \end{proof}
  \end{lem}

\subsection{Positivity and Legendrian polygons}

Given a generic, non-contractible Legendrian polygon $P$ with $2k$ vertices $p_1,\dots,p_{2k}$, we can associate to it a $k$-tuple of oriented isotropic flags as follows. Choose a set of representatives so that $P = P_-(\vec{v}_1,\dots,\vec{v}_{2k})$, and let $F_1 = (\Span(\vec{v}_1) \subset \Span(\vec{v}_1,\vec{v}_2))$, $F_2 = (\Span(\vec{v}_3) \subset \Span(\vec{v}_3,\vec{v}_4))$, \dots, $F_k = (\Span(\vec{v}_{2k-1}) \subset \Span(\vec{v}_{2k-1}, \vec{v}_{2k}))$. Since the representatives $\vec{v}_i$ are unique up to positive scalar multiplication or replacing all $\vec{v}_i$ with $-\vec{v}_i$, this procedure is well-defined on $P$. We denote the induced map by
\[\mathcal{F} : \mathscr{P}^-_{2k} \rightarrow \mathscr{F}^{(k)}.\]

\begin{rmk}
The configuration of oriented flags $\mathcal{F}(P)$ depends on the labeling of the vertices of $P$. If we re-label the vertices starting at any odd vertex $p_{2k+1}$ we obtain a cyclic permutation of the flags, but if we start at an even vertex or reverse the orientation we obtain different flags. In total, up to cyclic permutations and choice of cyclic ordering, there are four different (unlabeled) $k$-tuples of flags that could be associated to the Legendrian polygon $P$.
\end{rmk}

\begin{prop}
  Let $P = P_-(\vec{v}_1,\dots,\vec{v}_6)$ be a positive-transverse Legendrian polygon. Then,  the associated triple $\mathcal{F}(P) = (F_1,F_2,F_3)$ is positive.
  \begin{proof}
    We fix a basis $E$ in which the symplectic form is given by $\Omega$. Applying an element of $\PSp(\V,\w)$, we may assume that $F_1 = F_E$ and
    $F_3 = F_{\widehat{E}}$.
    
    Then, we may scale by positive scalars so that $\vec{v}_1 = \vec{e}_1$ and $\vec{v}_5 = \vec{e}_4$. Moreover, since $F_1^{(2)} = \Span(\vec{v}_1,\vec{v}_2)$, we have $\vec{v}_2 = a \vec{e}_1 + \vec{e}_2$ for some $a\in \bR$. Since $F_3^{(2)} = \Span(\vec{v}_5,\vec{v}_6)$ and $\w(\vec{v}_6,\vec{v}_1)=0$ we have $\vec{v}_6 = -\vec{e}_3$.
    
    Using the relations $\omega(\vec{v}_3,\vec{v}_5) > 0$ and $\omega(\vec{v}_4,\vec{v}_6) > 0$ we can normalize $\vec{v}_3$ and $\vec{v}_4$ and represent the flags $F_1,F_2,F_3$ by the following three matrices:
    \[\begin{pmatrix}
      1 & a\\
      0 & 1\\
      0 & 0\\
      0 & 0\\
    \end{pmatrix}, ~
    \begin{pmatrix}
      1 & 0\\
      b & 1\\
      a d & c\\
      d & e\\
    \end{pmatrix}, ~
    \begin{pmatrix} 
     0 & 0\\
     0 & 0\\
     0 & -1\\
     1 & 0\\
     \end{pmatrix},\]
     where $a,b,c,d,e \in \bR$, and $e - b c + a d = 0$. We must show that the matrix representing $F_2$ is totally positive. It suffices to show the positivity of $b$, $a d$, $d$, $b c - a d$, and $(a e - c)d$.
     
     But $b = \omega(\vec{v}_3,\vec{v}_6) > 0$, $a = \omega(\vec{v}_2,\vec{v}_5) > 0$, $d = \omega(\vec{v}_1,\vec{v}_3) > 0$, $e = b c - a d = \omega(\vec{v}_1, \vec{v}_4) > 0$, and $a e - c = \omega(\vec{v}_2,\vec{v}_4) > 0$, proving the claim.
  \end{proof}
\end{prop}

We get the following consequence for piecewise circular hexagons.

\begin{cor}
  A simple circular hexagon either has decreasing curvature or it has increasing curvature.
\end{cor}

For Legendrian polygons with more than $6$ vertices however, transversality does not suffice to guarantee the positivity of the associated tuple of flags.

\begin{figure}
    \centering
    \includegraphics[width=0.6\textwidth]{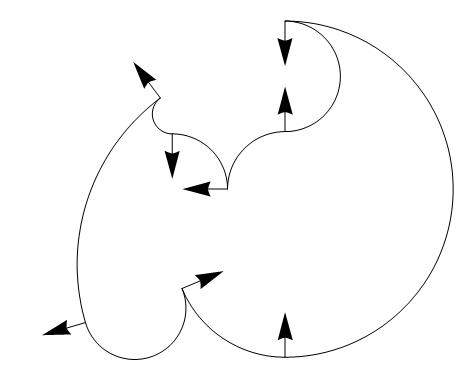}
    \caption{A simple circular octagon which is not positive.}
    \label{fig:transverseoctagon}
\end{figure}

An explicit example is the transverse Legendrian octagon with vertices
$(1,0,0,0)$, $(1,1,0,0)$, $(1,1,1,1)$, $(4,2,-1,1)$, $(4,8,-4,1)$, $(0,1,1,3)$, $(0,0,0,1)$, and $(0,0,-1,0)$ in a symplectic basis. We show its image in the space of oriented contact elements to $S^2$ in Figure \ref{fig:transverseoctagon}.

\begin{prop}\label{prop:polygonpositiveimpliesflagpositive}
  Suppose $P = P_{-}(\vec{v}_1,\dots,\vec{v}_{2N})$ is positive-transverse and has decreasing curvature. Then, the associated $N$-tuple $\mathcal{F}(P) = (F_1,\dots,F_N)$ is positive.
  \begin{proof}
    We first show that $F_1,F_2,F_k$ is positive for $k=3,\dots,N$.
    
    Choose a symplectic basis such that $\vec{v}_1 = \vec{e}_1$, $\vec{v}_2 = a \vec{e}_1 + \vec{e}_2$, $\vec{v}_{2k-1} = \vec{e}_4$ and $\vec{v}_{2k} = -\vec{e}_3 + b \vec{e}_4$. Then, normalizing each vector we may write $\vec{v}_1,\vec{v}_2,\vec{v}_3,\vec{v}_4,\vec{v}_{2k-1},\vec{v}_{2k}$ as the columns of the following matrix:
    \[\begin{pmatrix}
    1 & a & 1 & 1 &             0 & 0\\
    0 & 1 & c & e &             0 & 0\\
    0 & 0 & a d & f           & 0 & -1\\
    0 & 0 & d & d+c f - a d e & 1 & b\\
    \end{pmatrix}.\]
    So we must show that the triple of flags
    \begin{equation}\label{eqn:pospoly3flags}
    \begin{pmatrix}
     1 & 0\\
    0 & 1 \\
    0 & 0 \\
    0 & 0\\\end{pmatrix}
    ,~
    \begin{pmatrix}
    1 & 0 \\
    c & 1 \\
    a d & \frac{f-ad}{e-c}\\
    d & \frac{c f - a d e}{e-c}\end{pmatrix}
    ,~
    \begin{pmatrix}
    0 & 0\\
    0 & 0\\
    0 & -1 \\
    1 & 0\\\end{pmatrix}
    \end{equation}
    is positive, and for this it suffices to show that the second matrix is totally positive.
    
    From oriented transversality we obtain the inequalities:
    $d>0$, $a > 0$, $d + c f - a d e > 0$, 
    \begin{equation}\label{eqn:ineq0}
        a(d+c f - a d e)-f > 0.
    \end{equation}
    The decreasing curvature condition means that the Maslov index of the three Lagrangians (\ref{eqn:pospoly3flags}) is $1$. We find that the Maslov bilinear form is
    \[\begin{pmatrix}
     (a c-1) d & d (a e-1) \\
     d (a e-1) & d (a e-1)+(e-c) f \\
    \end{pmatrix}\]
    and so from Sylvester's positivity criterion we obtain
    \begin{equation}\label{eqn:ineq1}
        (a c  - 1)d > 0
    \end{equation}
    and
    \begin{equation}\label{eqn:ineqdet}
    d (e - c) (a(d + c f - a d e) - f) > 0.
    \end{equation}
    
    Inequality (\ref{eqn:ineq1}) together with $a>0$ and $d>0$ immediately implies $c>0$.
    
    Inequality (\ref{eqn:ineqdet}), together with $d>0$ and inequality (\ref{eqn:ineq0}) implies that $e-c>0$.
    
    Moreover, inequalities (\ref{eqn:ineq0}) and (\ref{eqn:ineq1}) as well as $d>0$ imply
    \[f > \frac{a d(a e - 1)}{a c - 1}.\]
    It follows that
    \begin{align*}
        c f - a d e &> \frac{c(ad)(a e - 1) - ade(ac-1)}{ac-1}\\
                   &= \frac{ad(e-c)}{a c - 1} > 0.
    \end{align*}
    Finally, using equation (\ref{eqn:ineq0})
    \[a d(d + c f - a d e) - d f > d f - d f = 0.\]
    This shows that the triple $F_1,F_2,F_k$ is positive.
    
    Since transversality and monotonicity are invariant under cyclic permutations, we might have chosen any odd index for $\vec{v}_1$ and so all triples of the form $F_i,F_{i+1},F_k$ with $i+1\neq k$ are positive. This implies, by Lemma \ref{lem:cyclicsetcharacterization}, that the $N$-tuple of flags is positive.
  \end{proof}
\end{prop}

\begin{figure}[h]
      \centering
      \includegraphics[width=.8\textwidth]{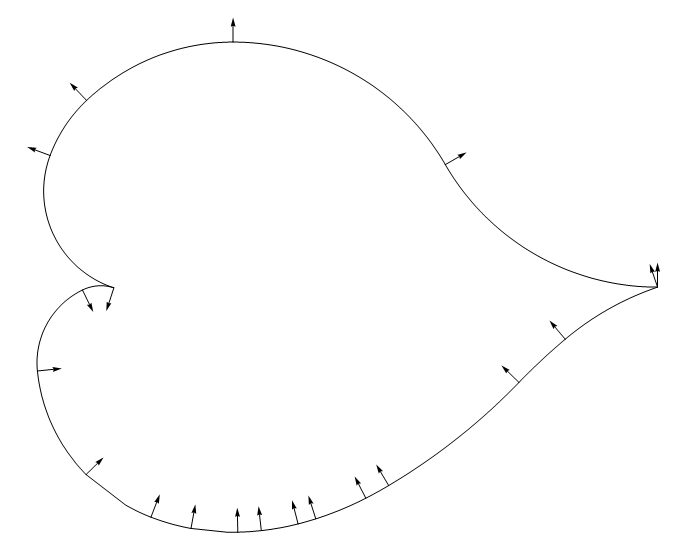}
      \caption{A simple circular $20$-gon with decreasing curvature. Oriented contact elements are indicated by arrows at the vertices.}
      \label{fig:convex20}
  \end{figure}

  We illustrate examples of simple circular $2k$-gons with decreasing curvature in Figures \ref{fig:threegraphs} and \ref{fig:convex20}.
  
  Next we will define an left inverse for $\mathcal{F}$. Given oriented subspaces $U,W\subset \V$, their intersection $U\cap W$ is endowed with the unique orientation satisfying: whenever $e_1,\dots,e_N$ is an oriented basis of $U+W$ such that $e_1,\dots,e_{k_2}$ is an oriented basis of $U$ and $e_{k_1},\dots,e_N$ is an oriented basis of $W$, then $e_{k_1},\dots,e_{k_2}$ is an oriented basis of $U\cap W$.
  
  Define a map $\mathcal{C} : \mathscr{F}_{>0}^{(k)} \rightarrow \mathscr{P}_{2k}$ by \begin{align*}
      \mathcal{C}(F_1,\dots,F_k) &= \\
      P_-(F_1^{(1)}, & ~F_1^{(2)}\cap F_2^{(3)}, F_2^{(1)}, \dots, F_{k-1}^{(2)}\cap F_{k}^{(3)}, F_k^{(1)}, F_k^{(2)}\cap F_1^{(3)}),
  \end{align*}
  
  where we slightly abuse notation since the image of $P_-$ does not depend on the choice of representative for each oriented $1$-dimensional subspace.
  
  \begin{example}
    The image by $\mathcal{C}$ of the parametrized triple of flags from Example \ref{ex:parametrization} is the Legendrian polygon $P_-(\vec{v}_1, \dots, \vec{v}_6)$, where
    \begin{align*}
    \vec{v}_1 = \begin{pmatrix} 1\\0\\0\\0\end{pmatrix},~
    \vec{v}_2 = \begin{pmatrix} y\\1\\0\\0\end{pmatrix},~
    \vec{v}_3 = \begin{pmatrix} 1\\x+y+\frac1y\\y\\1\end{pmatrix},\\
    \vec{v}_4 = \begin{pmatrix} 0\\1\\1\\x+\frac1y\end{pmatrix},~
    \vec{v}_5 = \begin{pmatrix} 0\\0\\0\\1\end{pmatrix},~
    \vec{v}_6 = \begin{pmatrix} 0\\0\\-1\\0\end{pmatrix}.\end{align*}
    It is simple to check that it is positive-transverse and has decreasing curvature.
  \end{example}
  
  \begin{prop}\label{prop:FChomeo}
    The map $\mathcal{F}$ is a homeomorphism onto its image.
    \begin{proof}
        Let $P = P_-(\vec{v}_1,\dots,\vec{v}_{2k})$ be a generic, noncontractible Legendrian polygon let $F_1,\dots,F_k = \mathcal{F}(P)$. Since $\w(\vec{v_{2j-1}},\vec{v_{2k-1}})\neq 0$ for all $j\neq k$, we have that $F_{j}^{(1)}$ is transverse to $F_{k}^{(3)}$ for all $j\neq k$. Moreover, since the edges spanned by $\vec{v}_{2j-1},\vec{v}_{2j}$ and $\vec{v}_{2k-1},\vec{v}_{2k}$ are linearly independent, $F_{j}^{(2)}$ is transverse to $F_{k}^{(2)}$ for all $j\neq k$. We conclude that the image of $\mathcal{F}$ is an $n$-tuple of pairwise transverse flags.
    
        The map $\mathcal{C}$, when restricted to the image $\mathcal{F}(\mathscr{P}^-_{2k})$, is an inverse of $\mathcal{F}$. Both maps are continuous since they are $\PSp(\V,\w)$-equivariant and the topologies on both moduli spaces are induced by that of the group, so the claim follows.
    \end{proof}
  \end{prop}

  \begin{prop}\label{prop:flagpositiveimpliespolygonpositive}
    Let $F_1,\dots,F_k$ be a positive $k$-tuple of oriented isotropic flags. Then, $\mathcal{C}(F_1,\dots,F_k)$ is a positive-transverse Legendrian $2k$-gon with decreasing curvature.
    \begin{proof}
      Positivity of the tuple $F_1,\dots,F_k$ means positivity of every cyclically ordered sub-triple. Consider a sub-triple of the form $F_i,F_{i+1},F_j$ with $j >i+1$ and apply an element of $\PSp(4,\bR)$ so that
      \[F_i = \begin{pmatrix} 1 & 0\\ 0 & 1\\0 & 0\\0 & 0\end{pmatrix}, \quad
      F_{i+1} = \begin{pmatrix} 1 & 0\\ a & 1\\ b & d \\ c & a d - b \end{pmatrix}, \quad 
      F_j = \begin{pmatrix} 0 & 0\\ 0 & 0\\ 0 & -1 \\ 1 & 0\end{pmatrix}.\]
      By positivity, $a,b,c>0$, $ad-b>0$, and  $-b^2+abd-cd>0$.
      Then, the vertices $\vec{v}_{2i-1},\vec{v}_{2i},\vec{v}_{2i+1},\vec{v}_{2i+2},\vec{v}_{2j-1},\vec{v}_{2j}$ of $\mathcal{C}(F_1,\dots,F_k)$ have coordinates of the form given in the columns of the following matrix:
      \[\begin{pmatrix}
      1 & \frac{b}{c} & 1 & k_1 & 0 & 0\\
      0 & 1           & a & 1 + k_1 a & 0 & 0\\
      0 & 0           & b & d + k_1 b & 0 & -1\\
      0 & 0           & c & ad-b + k_1 c & 1 & k_2\\
      \end{pmatrix}\]
      for some $k_1,k_2\in \bR$.
      
      We first show that $k_1,k_2\geq 0$.
      
      If $j=i+2$, then $k_1=0$. Otherwise, note that the flag $F_{i+2}$ has the form
      \[F_{i+2} = \begin{pmatrix} 1 & 0 & 0 & 0\\ x & 1 & 0 & 0\\ y & w & 1 & 0 \\ z & x w - y & x & 1 \end{pmatrix},\]
      and by positivity of the quadruple $F_0,F_{i+1},F_{i+2},F_j$ and Proposition \ref{prop:positivequadruple},
      \begin{align*}
      &\begin{pmatrix} 1 & 0 & 0 & 0\\ a & 1 & 0 & 0\\ b & d & 1 & 0 \\ c & a d - b & a & 1 \end{pmatrix}^{-1}
      \begin{pmatrix} 1 & 0 & 0 & 0\\ x & 1 & 0 & 0\\ y & w & 1 & 0 \\ z & x w - y & x & 1 \end{pmatrix} =\\
      &\begin{pmatrix} 1 & 0 & 0 & 0\\ x - a & 1 & 0 & 0\\ ad-b+y-dx & w-d & 1 & 0 \\ z - c - ay + bx & b-y+w(x-a) & x-a & 1 \end{pmatrix}\in U_-^{>0}.
      \end{align*}
      The fact that $\omega(\vec{v}_{2i+2},\vec{v}_{2i+3})=0$ implies that $k_1=\frac{ad-b+y-dx}{-ay+bx+z-c}$, which is positive since the numerator and denominator appear as entries in the first column of the above totally positive lower triangular matrix.
      
      Similarly, if $j=k$ and $i=1$ so that $F_j$  and $F_i$ are adjacent flags in the cyclic tuple, then $k_2=0$. Otherwise, the flag $F_{j+1}$ is of the form
      \[\begin{pmatrix}
      1 & 0 & 0 & 0\\
      -\alpha & -1 & 0 & 0\\
      \beta & \delta & 1 & 0\\
      -\gamma & -\alpha \delta + \beta & -\alpha & -1
      \end{pmatrix}\in KU_-^{>0}K,\]
      and $k_2 = \alpha > 0$.
      
      Now we compute the symplectic products: 
      \begin{align*}
          \omega(\vec{v}_{2i-1},\vec{v}_{2i+1}) = c &> 0; & \omega(\vec{v}_{2i-1},\vec{v}_{2i+2}) = ad - b + k c &> 0;\\
          \omega(\vec{v}_{2i-1},\vec{v}_{2i+2}) = ad - b + k c &> 0; & \omega(\vec{v}_{2i-1},\vec{v}_{2j-1}) = 1 &> 0;\\
          \omega(\vec{v}_{2i-1},\vec{v}_{2j}) = k_2 &\geq 0; & \omega(\vec{v}_{2i},\vec{v}_{2i+2}) = \frac{-b^2 + a b d - c d}{c} &> 0;\\
          \omega(\vec{v}_{2i},\vec{v}_{2j-1}) = \frac{b}{c} &> 0; & \omega(\vec{v}_{2i},\vec{v}_{2j}) = \frac{b k_2}{c} + 1 &> 0;\\
          \omega(\vec{v}_{2i+1},\vec{v}_{2j-1}) = 1 &> 0; & \omega(\vec{v}_{2i+1},\vec{v}_{2j}) = k_2 + a &> 0;\\
          \omega(\vec{v}_{2i+2},\vec{v}_{2j-1}) = k_1 &\geq 0; & \omega(\vec{v}_{2i+2},\vec{v}_{2j}) = k_1 k_2 + 1 + k_1 a &> 0.\\
      \end{align*}
      The two non-strict inequalities are strict as soon as they are not forced to be zero which occurs precisely when $F_j$ is adjacent to $F_i$ or $F_{i+1}$, in which case the vectors in the inequality represent adjacent vertices.
      
      Since $i,j$ are arbitrary, we conclude that $\omega(\vec{v}_i,\vec{v}_j)>0$ whenever $i+1 < j$ and $\mathcal{C}(F_1,\dots,F_k)$ is positive-transverse.
      
      Finally, we conclude from Lemma \ref{lem:PositiveMaslov} that the Maslov indices are $+1$, and so the Legendrian polygon has decreasing curvature.
    \end{proof}
  \end{prop}
  
  \begin{thm}
    The set of positive-transverse Legendrian $2k$-gons with decreasing curvature is a connected component of $\mathscr{P}_{2k}$ homeomorphic to a ball.
    \begin{proof}
      By  Proposition \ref{prop:FChomeo}, the map $\mathcal{F}:\mathscr{P}^-_{2k}\rightarrow \mathscr{F}^{(k)}$ is a homeomorphism onto its image. By Propositions \ref{prop:polygonpositiveimpliesflagpositive} and \ref{prop:flagpositiveimpliespolygonpositive}, the set of positive-transverse polygons with decreasing curvature is mapped onto the set of positive tuples of flags by $\mathcal{F}$. We conclude, by the Fock-Goncharov theorem (Theorem \ref{thm:fockgoncharov}, Proposition  \ref{prop:FockGoncharovComponent}), that the former is a connected component homeomorphic to a ball.
    \end{proof}
  \end{thm}
  
  Translating with the dictionary of Section \ref{sec:circles}, we find
  
  \begin{cor}
    The subspace of simple circular $2k$-gons with decreasing curvature is a connected component homeomorphic to a ball in the moduli space of generic circular $2k$-gons.
  \end{cor}

  \subsection{Flag curves}
  

  Let $\gamma : \bR \rightarrow \Proj_\bC(\V)$ be a smooth curve. A \emph{Legendrian lift} (or \emph{contact lift}) of $\gamma$ to $\Proj_\bR(V) \cong \mathbb{S}(T^*\Proj_\bC(\V))$ is a curve $\hat{\gamma} :\bR \rightarrow \Proj_\bR(V)$ such that $\pi \circ \hat{\gamma} = \gamma$. There are exactly two contact lifts of $\gamma$, given by a choice of transverse orientation.
  
  \begin{prop}
    Let $\gamma : \bR \rightarrow \Proj_\bC(\V)$ be a $C^2$ curve with non-vanishing curvature. Let $\hat{\gamma}$ be one of its two contact lifts to $\Proj_\bR(V)$. Then, the tangent lines to $\hat{\gamma}$ are Lagrangians which project to the osculating circles to $\gamma$.
    \begin{proof}
      For this proof we use coordinates as in Section \ref{sec:coordinates}. Assume without loss of generality that $\gamma$ is parametrized by arc length and write, in an affine chart,
      \[\gamma(s) = \begin{bmatrix}1\\z(s)\end{bmatrix}.\]
      Then, the contact lifts of $\gamma$ are
      \[\hat{\gamma}(s) = \sqrt{\pm\overline{T(s)}}\begin{pmatrix} 1\\
      z(s)\end{pmatrix},\]
      where $T(s) := z'(s)$ is the unit tangent vector to $z(s)$. Denote $T'(s)=k(s)N(s)$, where $N=iT$ is the unit normal vector and $k(s)$ the signed curvature function.
      We have
      \[\hat{\gamma}'(s) = \begin{pmatrix}
      \frac{\pm k(s)\overline{N(s)}}{2\sqrt{\overline{\pm T(s)}}}\\
      \frac{\pm k(s)\overline{N(s)}z(s)}{2\sqrt{\overline{\pm T(s)}}} + \sqrt{\overline{\pm T(s)}}T(s)
      \end{pmatrix},\]
      and a simple computation shows that $\w_\bC(\hat{\gamma}(s),\hat{\gamma}'(s)) = \pm i$, verifying that $\hat{\gamma}$ is Legendrian for $\w = \Re(\w_\bC)$.
      
      We also compute, choosing the positive sign in the choice of contact lift:
      \begin{align*}
          \pi(a \hat\gamma(s) + b\hat\gamma'(s)) &=
          \begin{bmatrix}
            2a\overline{T(s)} + b k(s)\overline{N(s)}\\
            z(s)(2a\overline{T(s)} + b k(s)\overline{N(s)}) + 2b
          \end{bmatrix}\\
          &=
          \begin{bmatrix}1\\ z(s) + \frac{2b}{2a\overline{T(s)} + bk(s)\overline{N(s)}}\end{bmatrix}\\
          &=\begin{bmatrix}1\\ z(s) + \frac{1}{k(s)}N(s)  - \frac{1}{k(s)}N(s) + \frac{4abT(s) + 2b^2 k(s)N(s)}{4a^2 + b^2k(s)^2}\end{bmatrix}\\
          &=\begin{bmatrix}1\\ z(s) + \frac{1}{k(s)}N(s) + \frac{4abT(s) + 2b^2 k(s)N(s)}{k(s)(4a^2 + b^2k(s)^2)}\end{bmatrix}\\
          &=\begin{bmatrix}1\\ z(s) + \frac{1}{k(s)}N(s) + \frac{1}{k(s)}\frac{(b k(s)-2a i)^2 N(s)}{(4a^2 + b^2k(s)^2)}\end{bmatrix}\\
      \end{align*}
      which, since $\left|\frac{(b k(s)-2a i)^2}{(4a^2 + b^2k(s)^2)}\right| = 1$, is a parametrization of the circle of radius $\frac{1}{k(s)}$ centered at $z(s) + \frac{1}{k(s)}N(s)$ in the affine patch. This circle is the osculating circle to the curve $z$ at $z(s)$. A similar computation with the negative sign yields the same circle (but corresponds to the opposite co-orientation).
    \end{proof}
  \end{prop}
  
  Let $\gamma : S^1 \rightarrow \Proj_\bR(\V)$ be a Legendrian curve. Motivated by the previous proposition, we call the projection to $\Proj_\bC(\V)$ of the tangent line at $\gamma(t)$ the \emph{osculating circle} to $\pi\circ \gamma$ at $\pi(\gamma(t))$.
  
\begin{figure}
      \centering
      \includegraphics[width=.7\textwidth]{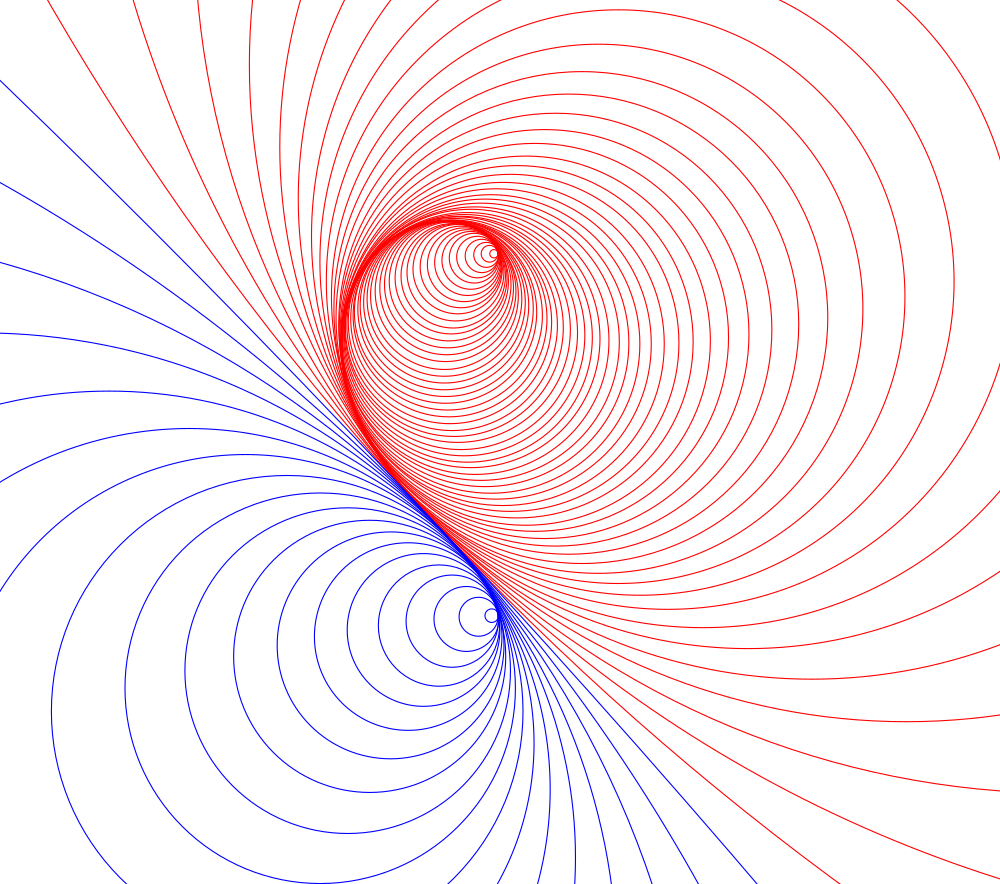}
      \caption{The osculating circles to half of the projected Veronese curve.}
      \label{fig:veronesecircleshalf}
  \end{figure}

  Accordingly, we say that $\gamma$ has \emph{decreasing curvature} if for every cyclically ordered $t_1,t_2,t_3\in S^1$, the triple of Lagrangians given by the tangent lines to $\gamma$ at $\gamma(t_i)$ have Maslov index $1$. See Figure \ref{fig:veronesecircleshalf} for an example of osculating circles to a decreasing curvature curve.
  
  Now consider a curve $\xi : S^1 \rightarrow \mathscr{F}$ in the space of isotropic flags.
  The curve $\xi$ is \emph{Frenet} if whenever $n_1+n_2+n_3 = p \le 3$, $0\le n_i\le 2$ we have
  \[\lim_{\substack{t_1,t_2,t_3 \text{ distinct}\\ t_i\rightarrow t}} \xi(t_1)^{(n_1)}\oplus\xi(t_2)^{(n_2)}\oplus\xi(t_3)^{(n_3)} = \xi(t)^{(p)}.\]
  
  \begin{prop}
    If $\xi$ is Frenet, the projection $\pi(\xi(t)^{(2)})$ is the osculating circle to the curve $\gamma = \pi\circ\gamma^{(1)}$ at $\gamma(t)$.
    \begin{proof}
      By definition of a Frenet flag curve, $\xi(t)^{(2)}$ is the tangent line to $\xi^{(1)}$ at $\xi(t)$, so its projection to $\Proj_\bC(\V)$ is the osculating circle to the projection $\pi\circ\xi^{(1)}$.
    \end{proof}
  \end{prop}
 
  A closed curve $\xi : S^1\rightarrow \mathscr{F}$ is \emph{positive} if for every cyclically ordered triple $t_1,t_2,t_3\in S^1$, the triple of flags $\xi(t_1),\xi(t_2),\xi(t_3)$ is positive. We obtain:
  
  \begin{thm}
    A positive Frenet curve in $\mathscr{F}$ is the tangent curve to the contact lift of an decreasing curvature simple closed curve in $\Proj_\bC(\V)$.
  \end{thm}
  
  By theorems of Labourie \cite[Theorem 1.4]{labourie} and Fock-Goncharov \cite[Theorems 1.14 and 1.15]{FG2006}, we know that the limit curve of a Hitchin representation is a positive Frenet curve.
  
  We conclude
  \begin{cor}
    The limit curve of a Hitchin representation in $\PSp(\V)$ is the contact lift of a simple closed curve in $\Proj_\bC(\V)$ with decreasing curvature.
  \end{cor}
  
\begin{figure}
      \centering
      \includegraphics[width=.5\textwidth]{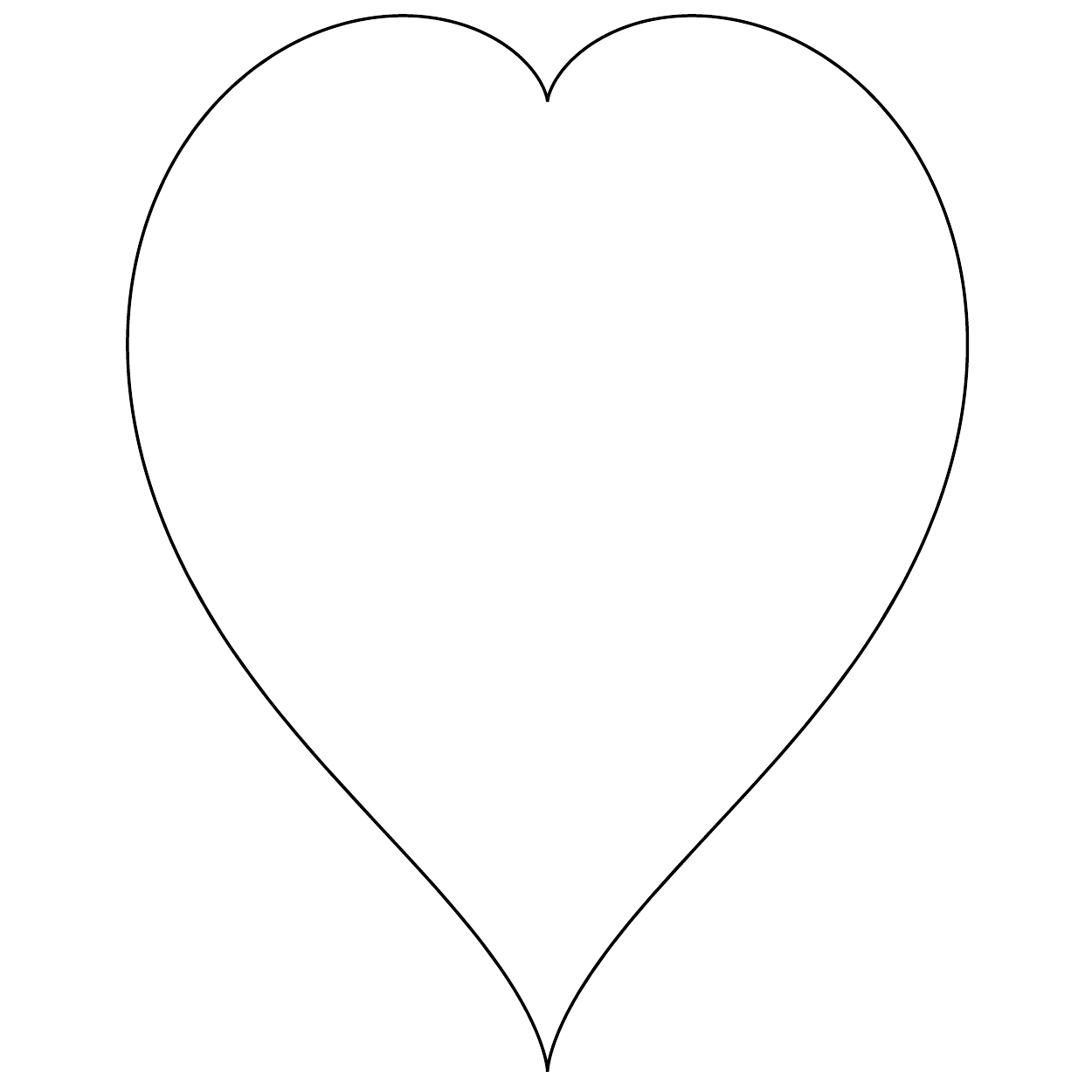}
      \caption{The projection to $\Proj_\bC(\V)$ of the Veronese curve.}
      \label{fig:veronese}
  \end{figure}

  For example, for a \emph{Fuchsian} representation, that is, one which factors through the irreducible representation $\PSL(2,\bR)\rightarrow \PSp(\V)$, the limit curve is the Veronese curve (or \emph{twisted cubic}) and its projection to $\Proj_\bC(\V)$ is depicted in Figure \ref{fig:veronese}.
  
\begin{figure}
      \centering
      \includegraphics[width=.7\textwidth]{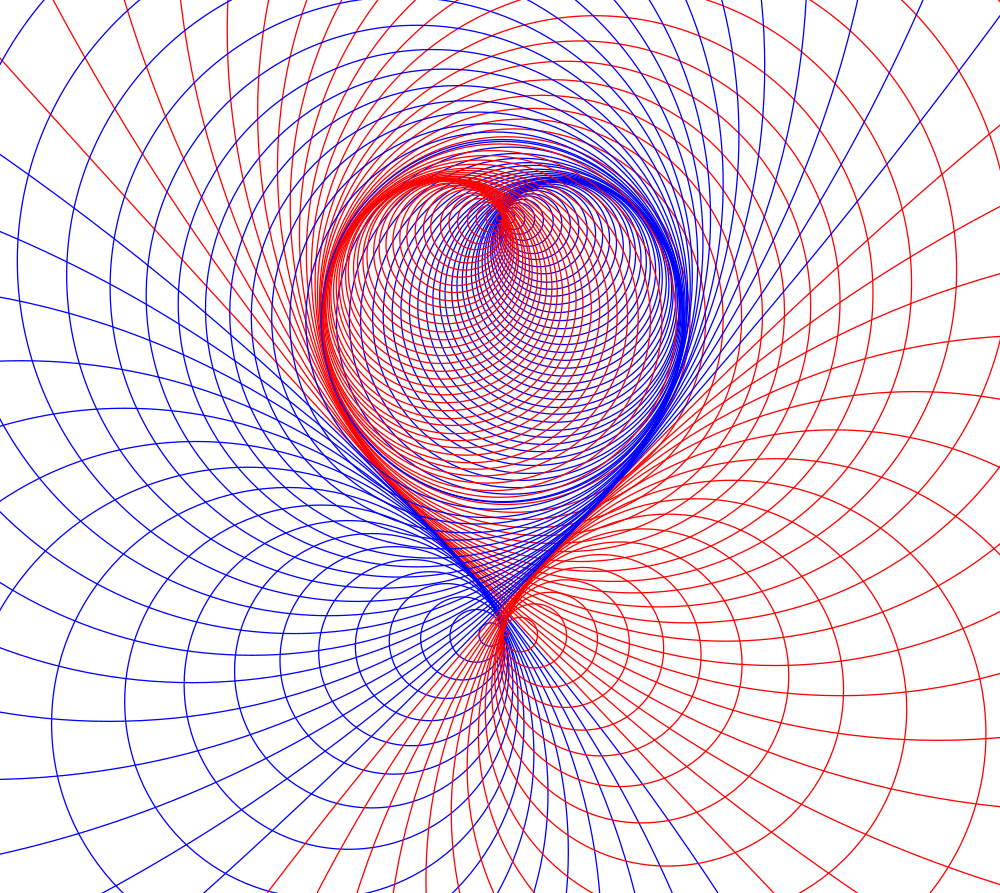}
      \caption{The osculating circles to the projected Veronese curve. Color represents orientation.}
      \label{fig:veronesecircles}
  \end{figure}

  We show the collection of osculating circles to the Veronese curve in Figure \ref{fig:veronesecircles}, and highlight the decreasing curvature by showing only half of them in Figure \ref{fig:veronesecircleshalf}.


\bibliographystyle{alpha}
\bibliography{bibliography.bib}
\end{document}